\documentclass{article}

\usepackage{amssymb}
\usepackage{amsmath}
\usepackage{theorem}
\usepackage{xypic}
\usepackage[all,v2]{xy}
\xyoption{2cell}
\UseAllTwocells
\xyoption{frame}

\allowdisplaybreaks[4]

\newtheorem{prop}{Proposition}[section]  
\newtheorem{lem}[prop]{Lemma}

\newtheorem{cor}[prop]{Corollary}
\newtheorem{them}[prop]{Theorem}

\newtheorem{conjecture}[prop]{Conjecture}

\theorembodyfont{\upshape}

\newtheorem{defn}[prop]{Definition}

\newtheorem{numrmk}[prop]{Remark}

\newtheorem{numex}[prop]{Example}
\newtheorem{numexamples}[prop]{Examples}

\newtheorem{note}{Note}

\newtheorem{example}{Example}

\newtheorem{rmk}{Remark}

\newenvironment{pf}{\begin{trivlist}\item[]{\sc Proof.}}%
            {\nolinebreak $\,\square$ \end{trivlist}}

\newcommand{\noprint}[1]{}

\renewcommand{\tilde}{\widetilde}
\newcommand{\dual}{^{\vee}}

\newcommand{\qed}{{\nolinebreak $\,\square$}}

\newcommand{\an}{{\mbox{\tiny an}}}

\newcommand{\stab}{{\mbox{\tiny stab}}}

\newcommand{\upst}{^{\ast}}

\newcommand{\lst}{_{\ast}}

\newcommand{\com}{^{\scriptscriptstyle\bullet}}

\newcommand{\argument}{{{\,\cdot\,}}}

\newcommand{\MM}{{\mathfrak M}}

\newcommand{\CC}{{\mathfrak C}}
\newcommand{\zz}{{\mathbb Z}}
\newcommand{\hh}{{\mathbb H}}

\newcommand{\aaa}{{\mathbb A}}

\newcommand{\pp}{{\mathbb P}}
\newcommand{\cc}{{\mathbb C}}

\newcommand{\Gm}{{{\mathbb G}_{\mbox{\tiny\rm m}}}}
\newcommand{\Ga}{{{\mathbb G}_{\mbox{\tiny\rm a}}}}

\newcommand{\tT}{{\cal T}}

\renewcommand{\O}{{\cal O}}

\newcommand{\cC}{{\cal C}}

\newcommand{\kK}{{\cal K}}

\newcommand{\hH}{{\cal H}}

\newcommand{\del}{\partial}
\newcommand{\resto}{{ \mid }}
\newcommand{\st}{\mathrel{\mid}}

\newcommand{\rk}{\mathop{\rm rk}}

\newcommand{\gr}{\mathop{\rm gr}\nolimits}

\newcommand{\Mor}{\mathop{\rm Mor}\nolimits}

\newcommand{\tr}{\mathop{\rm tr}\nolimits}

\newcommand{\sign}{\mathop{\rm sign}}

\newcommand{\spec}{\mathop{\rm Spec}\nolimits}

\newcommand{\Sym}{\mathop{\rm Sym}\nolimits}

\newcommand{\id}{\mathop{\rm id}}
\newcommand{\Hom}{\mathop{\rm Hom}\nolimits}

\newcommand{\Ext}{\mathop{\rm Ext}\nolimits}

\newcommand{\sheafhom}{\mathop{\rm {\mit{ \hH\! om}}}\nolimits}

\newcommand{\Isom}{\mathop{\rm Isom}\nolimits}

\newcommand{\Pic}{\mathop{\rm Pic}\nolimits}

\newcommand{\injectlim}%
           {\mathop{\lim\limits_{\textstyle\longrightarrow}}\limits}
\newcommand{\projectlim}%
           {\mathop{{\lim\limits_{\textstyle\longleftarrow}}}\limits}
\newcommand{\comp}{\mathbin{{\scriptstyle\circ}}}
\newcommand{\ol}{\overline}
\newcommand{\ul}{\underline}

\newcommand{\tensor}{\otimes}
\newcommand{\longiso}{\stackrel{\textstyle\sim}{\longrightarrow}}

\newcommand{\doublearrowstack}[2]%
           {{{{\scriptstyle#1}\atop{\textstyle\longrightarrow}}%
           \atop{{\textstyle\longrightarrow}\atop{\scriptstyle#2}}}}
\newcommand{\rightleftarrowstack}[2]%
           {{{{\scriptstyle#1}\atop{\textstyle\longrightarrow}}%
           \atop{{\textstyle\longleftarrow}\atop{\scriptstyle#2}}}}
\newcommand{\leftrightarrowstack}[2]%
           {{{{\scriptstyle#1}\atop{\textstyle\longleftarrow}}%
           \atop{{\textstyle\longrightarrow}\atop{\scriptstyle#2}}}}

\newcommand{\ses}[5]%
           {0\longrightarrow#1\stackrel{#2}{ \longrightarrow}#3\stackrel{#4}{
           \longrightarrow}#5\longrightarrow0}
\newcommand{\dt}[6]%
           {#1\stackrel{#2}{longrightarrow}#3 \stackrel{#4}{\longrightarrow}#5
           \stackrel{#6}{\longrightarrow} #1[1]}

\title{On the Cohomology of Stable Map Spaces}
\author{K. Behrend and A. O'Halloran}
\date{February 26, 2002}

\begin{document}
\sloppy
\maketitle

\begin{abstract}
We describe an approach to calculating the cohomology rings of stable
map spaces $\ol{M}_{0,0}(\pp^n,d)$.
\end{abstract}

\tableofcontents


\section*{Introduction}

Spaces of stable maps have enjoyed a lot of interest in recent years. They
were first introduced by Kontsevich in 1994 (see \cite{KM}), and have since
proven to be very useful in many contexts, especially in Quantum Cohomology
and Mirror Symmetry.  Stable map spaces are natural completions of spaces
of morphisms from algebraic curves to a fixed (non-singular complete)
variety $X$. They arise as natural generalizations of
the moduli spaces of stable curves discovered by Deligne and Mumford
\cite{DM}. 

Here we shall be mostly concerned with the spaces
$$\ol{M}_{0,0}(\pp^n,d)\,.$$ 
The generic member of $\ol{M}_{0,0}(\pp^n,d)$ (at least if $n\geq3$) is a
non-singular rational curve of degree $d$ in projective $n$-space and, in
fact, $\ol{M}_{0,0}(\pp^n,d)$ is a compactification of the space of all
such curves, and thus has dimension $dn+d+n-3$. The degenerations we allow
at the boundary are pairs $(C,f)$, where $C$ is a nodal curve of arithmetic
genus $0$ and $f$ is a morphism $f:C\to\pp^n$ of degree $d$, such that
every component of $C$ which is contracted to a point by $f$ has at least 3
nodes. This also 
gives the correct picture for $n=1,2$.  For example,
$\ol{M}_{0,0}(\pp^1,d)$ compactifies the space of degree $d$ ramified
covers of genus zero of the projective line. We also remark that
$\ol{M}_{0,0}(\pp^n,1)$ is simply the Grassmannian $G(2,n+1)$ of lines in
$\pp^n$.

The true beauty of the spaces $\ol{M}_{0,0}(\pp^n,d)$ only becomes
apparent if we consider them as {\em stacks}.  In fact, the algebraic
stacks underlying the various $\ol{M}_{0,0}(\pp^n,d)$ are smooth and admit
universal families. These are properties that the spaces
$\ol{M}_{0,0}(\pp^n,d)$ generally lack.  We shall always work with these
stacks and thus use the notation $\ol{M}_{0,0}(\pp^n,d)$ for the {\em stack
of stable maps }of degree $d$ to $\pp^n$ (of genus $0$ without marked
points).

Our goal is to compute the cohomology ring of $\ol{M}_{0,0}(\pp^n,d)$, or
at least, to outline an approach by which this might be
achieved. (The only previous result in this direction is the computation of
the Betti numbers of $\ol{M}_{0,\nu}(\pp^n,d)$ due
to Getzler and Pandharipande \cite{getzler}.) Our
method is inspired by the utility of $\cc^\ast$-actions for studying
integrals over stable map spaces, but there is an additional ingredient: a
vector field (which is compatible with the $\cc^\ast$-action).

This method is due to Akildiz and Carrell \cite{carrell} and can be
summarized as follows.  Let $X$ be a non-singular projective variety
over $\cc$ with a
$\cc^\ast$-action and suppose that $V$ is a vector field on $X$, satisfying
${^\lambda}V=\lambda V$, for all $\lambda\in\cc^\ast$ (we say that $V$ is
{\em equivariant\/}).  If $V$ has exactly one fixed point and $Z$ is
the scheme-theoretic fixed locus of $V$ (so $Z$ is a one-point non-reduced
scheme), then we have
$$H^\ast(X,\cc)=\Gamma(Z,\O_Z)\,.$$
(See Examples~\ref{stupex} and~\ref{morestupex}, below, where this is
worked out for the special case of $X=\pp^n$. See also
Remark~\ref{final}, for the case of the Grassmannian of lines in $\pp^n$.)

Our use of the method of Akildiz-Carrell is novel in two aspects: we apply
it to stacks, but more significantly, the fixed locus
$Z$ of the vector field $V$ has positive dimension.  Thus we have to
replace the ring of global sections $\Gamma(Z,\O_Z)$ by the
hypercohomology ring $\hh^0(X,K_V\com)$, where $K_V\com$ is the Koszul
complex defined by the vector field $V$.

There is one important difference between $\Gamma(Z,\O_Z)$ and
$\hh^0(X,K_V\com)$. The ring $\Gamma(Z,\O_Z)$ can be computed entirely on
the fixed locus $Z$, whereas $\hh^0(X,K_V\com)$ depends on an open
neighborhood of $Z$ in $X$.  Thus, in the case of positive-dimensional
fixed locus $Z$, the localizing power of the method is much weaker. 

The method is saved by a somewhat surprising phenomenon.
We discovered that we can restrict our attention entirely to a certain
open subset $U$ of $X$, even though $U$ does {\em not} cover the fixed
locus $Z$ completely. This open subset $U$ is the `big cell'
associated by Bia\l ynicki-Birula \cite{BB} to the $\cc^\ast$-action
on $X$.

In the cases we consider here, it turns out that the canonical map
$$H^\ast(X,\cc)=\hh^0(X,K_V\com)\longrightarrow\hh^0(U,K_V\com)=
\Gamma(U,\O_Z)$$ 
is, though not injective, injective in all {\em relevant } degrees.
This means injective in all degrees that contain a generator or a
relation. 

One of our main results is an explicit description of the big Bia\l
ynicki-Birula cell of $\ol{M}_{0,0}(\pp^n,d)$ as a vector bundle over
$\ol{M}_{0,d}$ (modulo an action of the symmetric group $S_d$). Here
$\ol{M}_{0,d}$ is the space of stable curves of genus zero with $d$ marked
points, which is comparatively well understood. 

The case $d=3$ is particularly simple and we focus on it in the latter
part of the paper. If $d=3$ then $\ol{M}_{0,d}$ is just a point and so the
Bia\l ynicki-Birula cell is simply an affine space $\aaa^{4n}$ modulo an
action of $S_3$.  We succeed in writing down the vector field $V$ in
canonical coordinates on $\aaa^{4n}$.  This leads at least to a conjectural
description of the cohomology ring of $\ol{M}_{0,0}(\pp^n,3)$.  The truth
of this conjecture depends only on a certain purely algebraic statement,
which we verified using Macaulay~2 \cite{M2} for $n\leq5$. 

More interesting than the case of finite $n$ is the case of the limit as
$n$ approaches $\infty$.  The cohomology ring of $\ol{M}_{0,0}(\pp^n,d)$
stabilizes as $n$ increases, so we can define a ring which we call
the cohomology ring of $\ol{M}_{0,0}(\pp^\infty,d)$, even though this
latter stack does not make sense. 

We succeed in describing the cohomology ring of
$\ol{M}_{0,0}(\pp^\infty,3)$ completely using generators and
relations (Theorem~\ref{themi}).  This is the main result of the
paper. It says
$$H^\ast\big(\ol{M}_{0,0}(\pp^\infty,3),\cc\big)=
\cc[b,\sigma_1,\rho,\sigma_2,\tau,\sigma_3]/
\big((\tau^2-\rho\sigma_2),\tau\sigma_3,\rho\sigma_3\big)\,.$$
The generators can be expressed in terms
of Chern classes of certain canonical vector bundles on
$\ol{M}_{0,0}(\pp^\infty,3)$. The degrees of $b$ and $\sigma_1$ are 1,
the degrees of $\rho$, $\sigma_2$ and $\tau$ are 2 and the degree of
$\sigma_3$ is 3, using algebraic degrees (where the first Chern class
has degree 1).  Thus the degrees of the relations are 4, 5 and 5,
respectively. 

Thus, as a ring,
$H_{DR}\big(\ol{M}_{0,0}(\pp^\infty,3)\big)$ is reduced, of pure dimension
4 and has 
two irreducible components, $\cc[b,\sigma_1,\sigma_2,\sigma_3]$ and
$\cc[b,\sigma_1,\rho,\sigma_2,\tau]/(\tau^2-\rho\sigma_2)$, intersecting
transversally along $\cc[b,\sigma_1,\sigma_2]$. 

The case $d=2$ is special, as $\ol{M}_{0,2}$ is not defined. It is
much easier than the case $d=3$ and we have complete results.

We briefly outline the structure of the paper.

In Section~\ref{outline} we describe the theory of equivariant vector
fields and their relation to de Rham cohomology. We verify that the results
of Akildiz and Carrell which we require hold for stacks.  We improve
on existing treatments of Chern classes by proving that the
Carrell-Lieberman class \cite{CL} is {\em homogeneous }(see
Section~\ref{Chern}). Hence the Carrell-Lieberman characteristic classes
(and not just their leading terms) are equal to the corresponding Chern
classes.

Section~\ref{prelims} assembles a few facts about stable map stacks which
we require later. We observe that if $E$ is a convex vector bundle on the
variety $X$, then the stack of stable maps to $E$ is a vector bundle over
the stack of stable maps to $X$. We prove that the cohomology of
$\ol{M}_{0,\nu}(\pp^n,d)$ stabilized as $n$ 
increases and we define the cohomology ring of
$\ol{M}_{0,\nu}(\pp^\infty,d)$. Moreover,this cohomology
ring maps surjectively onto the cohomology ring of
$\ol{M}_{0,\nu}(\pp^n,d)$, for every finite $n$.

In
Section~\ref{SecPar} we describe the big Bia\l
ynicki-Birula cell of $\ol{M}_{0,0}(\pp^n,d)$ as a vector bundle over
$[\ol{M}_{0,d}/S_d]$. The most significant case is $n=1$. Here the big
Bia\l ynicki-Birula cell consists of all stable maps which are
unramified over $\infty\in\pp^1$.  In the general case, it consists of
all stable maps which avoid the codimension 2 plane
$\langle0,0,\ast,\ldots,\ast\rangle$ and intersect the hyperplane
$\langle 0,\ast,\ldots,\ast\rangle$ transversally $d$ times.

By changing the $\cc^\ast$-action on $\pp^n$, we can cover all of
$\ol{M}_{0,0}(\pp^n,d)$ with vector bundles over $\ol{M}_{0,d}$, and so our
results lead, at least in principle, to an explicit description of the stable
map stacks $\ol{M}_{0,0}(\pp^n,d)$ in terms of stable curve spaces
$\ol{M}_{0,d}$. 

Section~\ref{thevectorfield} starts with a recipe to calculate our vector
field on $\ol{M}_{0,0}(\pp^n,d)$.  The key result is that the derivative of
the universal map $f:\tilde{C}\to\pp^n$ induces an isomorphism
$\Gamma(\tilde{C},\tT_{\tilde{C}})\to
\Gamma(\tilde{C},f\upst\tT_{\pp^n})$. The remainder of the section contains
the calculations for the cases $d=2$  and $d=3$.

\subsection*{Notation and Conventions}

Throughout the paper we will work over the ground field $\cc$ of
complex numbers.

All of our algebraic stacks will be of Deligne-Mumford type.  This
means that the diagonal $X\to X\times X$ is unramified.
Deligne-Mumford stacks $X$ admit \'etale presentations $U\to X$, where
$U$ is a scheme.  We denote the stack quotient associated to a
$G$-variety $X$ by $[X/G]$.

Whenever we consider sheaves on a Deligne-Mumford stack $X$, it is
understood that these are sheaves on the small
\'etale site of $X$. Objects of this \'etale site are thus \'etale
morphisms $U\to X$, where $U$ is a scheme.  The topology on this site
is defined in the same way as for the \'etale site of a scheme. Any
vector bundle $E\to X$ defines a sheaf of local sections on the
\'etale site of $X$, which we often identify with $E$.

Any cohomology group of a sheaf on $X$ is understood to be the
cohomology of the \'etale site with values in the given sheaf, unless
mentioned otherwise.

\subsubsection*{Stable maps}

For an algebraic variety $X$ (not necessarily proper), we denote by
$H_2(X)^+$ the semigroup (with $0$) of group homomorphisms $\Pic(X)\to\zz$,
which take non-negative values on ample line bundles. This semigroup is a
convenient set of labels for the class of a stable map.  Given a stable map
$(C,x,f)$ to $X$, it is of class $\beta\in H_2(X)^+$ if $\deg(f\upst
L)=\beta(L)$, for all $L\in\Pic(X)$. If $X=\pp^n$, we identify $H_2(X)^+$
with $\zz_{\geq0}$.

We denote by $\ol{M}_{g,n}(X,\beta)$ the stack of stable maps of class
$\beta$ from $n$-marked genus $g$ curves to $X$.

\subsubsection*{Group actions}

If an algebraic group $G$ acts on a smooth scheme or a smooth algebraic
stack $X$, we will denote this action on the right. Assume given a right
action of $G$ on $X$ and a lift of this action to a vector bundle $E$ over
$X$. Then for any $g\in G$ and any local section $e\in\Gamma(Ug,E)$, we
denote by ${^g}e$ the section of $E$ over $U$ given by the formula
\begin{equation}\label{funact}
({^g}e)(x)=e(xg)g^{-1}\,.
\end{equation}
In other words, ${^g}e\in\Gamma(U,E)$ is defined to make the diagram
$$\xymatrix{
E\rto^g& E\\
U\rto^g\uto^{{^g}e}& Ug\uto_{e}}$$
commute.
In particular, (\ref{funact}) defines a left representation of $G$
on the 
$\cc$-vector 
space $\Gamma(X,E)$ of global sections of $E$.

Any action of $G$ on $X$ lifts naturally
to the vector bundles $\O_X$, $\tT_X$ and $\Omega_X$, and so we get induced
representations of $G$ on functions, vector fields and differential forms on
$X$. Explicitly, if $f\in\Gamma(X,\O_X)$ is a regular function on $X$,
then ${^g}f=g\upst f$, or $({^g}f)(x)=f(xg)$, for all $x\in X$. If
$V\in\Gamma(X,\tT_X)$ is a 
vector field on $X$, then ${^g}V$ is characterized by the formula 
$$Dg(x)\big(({^g}V)(x)\big)=V(xg)\,$$
or more briefly by 
$$(Dg)({^g}V)=g\upst V\,.$$
(Here $Dg:\tT_X\to g\upst\tT_X$ is the derivative of $g:X\to X$.)
If $\omega\in\Gamma(X,\Omega)$ is a differential form on $X$, then
${^g}\omega$ is given by 
$${^g}\omega=dg(g\upst\omega)\,,$$
or 
$$({^g}\omega)(x)=dg(x)\big(\omega(xg)\big)\,,$$
where $dg:g\upst\Omega_X\to\Omega_X$ denotes the natural pullback
homomorphism. Note that for all $g\in G$ we have
\begin{equation}\label{scal}
{^g}\langle\omega,V\rangle=\langle{^g}\omega,{^g}V\rangle\,.
\end{equation}

If $F$ is another vector bundle over $X$
to which the $G$-action has been lifted, then we get induced $G$-actions
also on the vector bundles $\sheafhom(E,F)$ and $E\otimes F$, given by the
formulas $(\phi g)(eg)=\phi(e)g$ and $(e\otimes f)g=eg\otimes fg$. On
global sections, this gives rise to $G$-representations by the formulas
$$({^g}\phi)({^g}e)={^g}\big(\phi(e)\big)$$
and 
$${^g}(e\otimes f)={^g}e \otimes{^g}f\,.$$

Finally, given a $\cc$-linear sheaf homomorphism $\Phi:E\to F$, we define
${^g}\Phi$ by
$$({^g}\Phi)({^g}e)={^g}\big(\Phi(e)\big)\,,$$
for any local section $e$ of $E$.  This generalizes the definitions above
if $\Phi$ is $\O_X$-linear and gives rise to a $G$-representation
on the space of all $\cc$-linear sheaf homomorphisms from $E$ to $F$. For
example, the universal derivation $d:\O_X\to\Omega_X$ satisfies ${^g}d=d$,
for all $g\in G$ (this follows easily from~(\ref{scal})). 
Note that, because of this, ${^g}\nabla$, for a connection $\nabla$ on $E$,
is again a connection on $E$.

Whenever any kind of object $A$ satisfies an equation ${^g}A=A$, for all
$g\in G$, then we refer to $A$ as {\em $G$-invariant}. 

\subsubsection*{The case of $\Gm$}\label{gmcase}

In many cases our group $G$ will be equal to the multiplicative group
$\Gm$. If we are given a right $\Gm$-action on a vector bundle $E$, lifting
a right action of $\Gm$ on $X$, then we refer to this action as the {\em
geometric }action, to distinguish it from the action of $\Gm$ on $E$ by
scalar multiplication on the fibers, which we shall call the {\em linear
}action. Of course the geometric and the linear action commute with each
other. 

In this case the (geometric) $\Gm$-action gives rise to a
$\cc^\ast$-representation on all the above mentioned vector spaces, i.e.,
it makes them into graded $\cc$-vector spaces.  The homogeneous elements
$A$ of degree $i$ satisfy
$${^\lambda}A=\lambda^i A\,,$$
for all $\lambda\in\cc^\ast$. Of particular interest to us are elements of
degree one; we shall call them {\em $\Gm$-equivariant}, or just {\em
equivariant}, since there is rarely any fear of confusion.
In particular, this gives rise to the notion of {\em equivariant vector
field}.


\section{Outline of the Method}\label{outline}

Let $X$ be a smooth and proper Deligne-Mumford stack. Assume that 
$$\text{$H^{p}(X,\Omega^q)=0$, for all $p\not=q$,}$$
where $\Omega=\Omega_X$ is the sheaf of K\"ahler differentials on $X$.
(See Remark~\ref{pqnz}.)

We will be interested in the graded ring
$$H_{DR}(X)=\bigoplus_p H^p(X,\Omega^p)\,\,.$$
The notation $H_{DR}(X)$ is justified, because under our assumption
$$\bigoplus_p H^p(X,\Omega^p)=\hh\big(X,(\Omega\com,d)\big)\,\,,$$
where $(\Omega\com,d)$ is the algebraic de Rham complex of $X$.
We use the algebraic grading, i.e., we consider $H^p(X,\Omega^p)$ to
have degree $p$.

\begin{numrmk}\label{allthecoh}
Let ${X_\an}$ be the (small) analytic site of $X$.
By a theorem of Grothendieck (see \cite{GrothDeRham}) we have
$$\hh\big(U,(\Omega\com,d)\big) = H(U_\an,\cc),$$
for every smooth variety $U$. Choosing an \'etale  presentation $U\to
X$, where $U$ is a smooth variety, we obtain $E_2$-spectral sequences
$$\hh^q\big(U_p,(\Omega\com,d)\big)
\Rightarrow\hh^{p+q}\big(X,(\Omega\com,d)\big)$$
and
$$H^q({U_p}_\an,\cc) 
\Rightarrow H^{p+q}(X_\an,\cc)\,,$$
where $U_p$ is the $(p+1)$-fold fibered product of $U$ with itself
over $X$. 

Thus we conclude that 
$$\hh\big(X,(\Omega\com,d)\big)=H({X_\an},\cc)\,\,.$$
Letting $\ol{X}$ be the coarse moduli space of $X$, we have
$$H({X_\an},\cc)=H({\ol{X}_\an},\cc)\,\,,$$
essentially because group cohomology of any finite group with values in
$\cc$ vanishes.
Thus we have
$$H_{DR}(X)=H({X_\an},\cc)=H({\ol{X}_\an},\cc)\,\,,$$ and so we can
also interpret our ring $H_{DR}(X)$ as the usual (singular)
cohomology ring of the topological space underlying the variety
$\ol{X}$ over $\cc$. 

Note that for $H(X_\an,\cc)$ we also use the algebraic grading, i.e.,
we consider $H^{2p}(X_\an,\cc)$ to have degree $p$.
\end{numrmk}

\subsection{Equivariant vector fields}

Now assume that we are given a (right) $\Gm$-action on $X$.
Recall that we defined a vector field $V$ on $X$ to be
$\Gm$-equivariant if ${^\lambda}V=\lambda V$, for all
$\lambda\in\cc^\ast$. Note that $V$ is equivariant if and only if 
the diagram 
$$\xymatrix{
\lambda^\ast\Omega_X\rto^{\lambda^\ast V}\dto_{\lambda d\lambda} &
\lambda^\ast\O_X \dto^{=} \\
\Omega_X\rto^{V} & \O_X}$$
commutes, for all $\lambda\in\cc^\ast$.

\begin{numrmk}\label{pair}
Let $G$ be the semidirect product of the multiplicative group and the
additive group, where the action of $\Gm$ by conjugation on $\Ga$ is
given by scalar multiplication: $^\lambda a=\lambda a$.  We can
identify $G$ with the group of $2\times2$ invertible matrices of the
form
$$\left(\begin{array}{cc}
1&0\\a&\lambda
\end{array}\right) \begin{array}{c}\phantom{.}\\ . \end{array} $$
Suppose we are given a right action of $G$ on $X$. Then restricting to
the multiplicative subgroup of $G$ gives us a (right) action of $\Gm$
on $X$.  Taking the derivative of the $\Ga$-action on $X$ defines a
vector field $V$ on $X$.  More precisely, let 
\begin{align*}
\phi_x:\aaa^1 & \longrightarrow X \\*
a & \longmapsto x a
\end{align*}
denote the orbit map of $x\in X$ for the  action of $\Ga$. Then
$V(x)=D\phi_x(0)$, where we identify the linear map 
$$D\phi_x(0):\tT_{\aaa^1}(0)\longrightarrow \tT_X(x)$$
with the image of the canonical generator $1\in \tT_{\aaa^1}(0)$.

Taking the derivative at $0$ of the commutative diagram
$$\xymatrix{
\aaa^1\dto_\lambda\rto^{\phi_{x\lambda}} & X \\
\aaa^1\rto^{\phi_x} & X \uto_{\lambda} }$$
proves that $D\phi_{x\lambda}(0)=D\lambda(x)D \phi_x(0)\lambda$, and
hence that ${^\lambda} V(x)=\lambda V(x)$.  Thus $V$ is
$\Gm$-equivariant.  This is the most common source of $\Gm$-equivariant
vector fields.
\end{numrmk}

\begin{numex}\label{stupex}
Let $D_\lambda$ denote the diagonal $(n+1)\times(n+1)$-matrix with
entries $(1,\lambda,\ldots,\lambda^n)$ along the diagonal. Let $N$
denote the nilpotent $(n+1)\times(n+1)$-matrix with ones along the
sub-diagonal and zeros elsewhere. Then
$$\left(\begin{array}{cc}
1&0\\a&\lambda
\end{array}\right) \longmapsto e^{aN} D_\lambda$$
defines a group homomorphism $G\to GL(n+1)$. Via this homomorphism, we
define a right action of $G$ on $\pp^n$, by acting in the natural way
on homogeneous coordinate vectors of $\pp^n$, which we consider to be
row vectors of length $n+1$.

As in Remark \ref{pair}, we get an induced $\Gm$-action and an
equivariant vector field $W$ on $\pp^n$. In standard homogeneous
coordinates on $\pp^n$ this vector field is given by 
$$W=\sum_{i=1}^n x_i \frac{\del}{\del  x_{i-1}}\,.$$

The zero locus of this vector field consists of one point, namely the
origin in the standard affine space $\aaa^n\subset\pp^n$ defined by
$x_0=1$. In affine coordinates $s_i=\frac{x_i}{x_0}$ this vector field
is given by 
$$W=-s_1s_n\frac{\del}{\del s_n} + \sum_{i=1}^{n-1}(s_{i+1}-s_1s_i)
\frac{\del}{\del s_i}\,.$$
(Recall the relation $\sum_{i=0}^nx_i\frac{\del}{\del x_i}=0$.)
\end{numex}

\subsubsection{The Koszul complex}

Let $X$ be as above, endowed with a $\Gm$-action. Let
$V:\Omega_X\to\O_X$ be an equivariant vector field on
$X$.  We can associate to $V$ the Koszul complex
$$\Omega^N_X\stackrel{\iota (V)}{\longrightarrow}
\ldots\stackrel{\iota(V)}{\longrightarrow}
\Omega_X^2\stackrel{\iota(V)}{\longrightarrow}
\Omega_X\stackrel{V}{\longrightarrow}\O_X\,,$$ where $N$ denotes the
dimension of $X$, hence the rank of $\Omega_X$, and $\iota(V)$ denotes
contraction with $V$. (If we pull back this Koszul complex to an
\'etale and affine $X$-scheme over which we can trivialize $\Omega_X$,
then the vector field $V$ is given by $N$ regular functions and the
above complex is the usual Koszul complex associated to this sequence
of regular functions.)

We set $K^p=\Omega_X^{-p}$, for $p\in\zz$, and denote the above
Koszul complex by 
$$K_V\com=\big(K\com,\iota(V)\big)\,.$$
Note that $K_V\com$ is a sheaf of differential graded commutative (with
unit) $\O_X$-algebras on $X$.

We are interested in the hypercohomology $\hh(X,K_V\com)$, and its relation
to $H_{DR}(X)$.

\begin{lem}
For every $p\not=0$ we have $\hh^p(X,K_V\com)=0$. Moreover,
$\hh^0(X,K_V\com)$ is a commutative filtered $\cc$-algebra, and for the
associated graded algebra we have
$$\gr\hh^0(X,K_V\com)=H_{DR}(X)\,.$$
\end{lem}
\begin{pf}
This follows immediately from the standard $E_1$ spectral sequence of
hypercohomology and our assumption on the vanishing of off-diagonal
Hodge groups. 
\end{pf}

For $\lambda\in\cc^*$ we define a homomorphism 
$\psi_\lambda:\lambda^\ast K_V\com\to K_V\com$ by
\begin{equation}\label{psilam}
\vcenter{\xymatrix{
\ldots\rto&
\lambda^\ast\Omega_X^p
\rto^{\lambda^\ast\iota(V)}\dto_{\lambda^{p}d\lambda} &
\ldots\rto &
\lambda^\ast\Omega_X\rto^{\lambda^\ast V}\dto_{\lambda d\lambda} &  
\lambda^\ast\O_X \dto^{=} \\
\ldots\rto&
\Omega_X^p\rto^{\iota(V)} &\ldots\rto &\Omega_X\rto^{V} & \O_X}}
\end{equation}
Here $d\lambda:\lambda^\ast\Omega^p_X\to\Omega^p_X$ is the
homomorphism induced on the exterior power by the derivative
$d\lambda:\lambda^\ast\Omega_X\to\Omega_X$. Note that $\psi_\lambda$ is a
homomorphism of differential graded algebras.

We can use $\psi_\lambda$ to define a representation of $\cc^*$ on the
hypercohomology $\hh^0(X,K_V\com)$. We simply associate to
$\lambda\in\cc^*$ the automorphism
$$\hh^0(X,K_V\com)\stackrel{\lambda^\ast}{\longrightarrow}
\hh^0(X,\lambda^\ast K_V\com) \stackrel{\psi_\lambda}{\longrightarrow}
\hh^0(X,K_V\com)\,.$$ 
Thus $\hh^0(X,K_V\com)$ becomes a graded $\cc$-algebra.

\begin{prop}[Akildiz-Carrell \cite{carrell}]\label{AC}
There is a canonical isomorphism of graded $\cc$-algebras
$$\hh^0(X,K_V\com)=H_{DR}(X)\,.$$
\end{prop}
\begin{pf}
Note that the grading induced by the $\cc^*$-representation on
$\hh^0(X,K_V\com)$ is compatible with the filtration induced by the
$E_1$-spectral sequence. Thus we get an induced $\cc^*$-representation on
the associated graded algebra. One shows that via this induced
representation, $\cc^*$ acts on $H^p(X,\Omega^p)$ through the character
$\lambda\mapsto\lambda^p$. Then the proof is finished, in view of the lemma
from linear algebra which we state below.

For the convenience of the reader, we recall the proof that $\cc^*$ acts on
$H^p(X,\Omega^p)$ through $\lambda\mapsto\lambda^p$. One simply factors
$\psi_\lambda$ as $\psi_\lambda=\phi_\lambda\comp d\lambda$, where
$\phi_\lambda$ is multiplication by $\lambda^{p}$ on
$\Omega^p$. Thus the action of $\lambda\in\cc^*$ on
$H^p(X,\Omega^p)$ factors as $\psi_\lambda\comp\lambda^\ast=
\phi_\lambda\comp(d\lambda\comp\lambda^\ast)$. Now
$d\lambda\comp\lambda^\ast$ is the homomorphism on de Rham cohomology
induced by the morphism $\lambda:X\to X$, which is the identity.  On the
other hand, $\phi_\lambda$ obviously induces multiplication by $\lambda^p$
on $H^p(X,\Omega^p)$. 
\end{pf}

\begin{lem}
Let $H$ be a commutative filtered $\cc$-algebra with a $\cc^*$-representation,
respecting the filtered algebra structure. Denote the filtration by 
$$\ldots\subset F_{i-1}\subset F_i\subset\ldots$$
Suppose that $\cc^*$ acts on $F_i/F_{i-1}$ through
$\lambda\mapsto\lambda^i$. Let $H=\bigoplus_j H_j$
be the grading induced by the $\cc^*$-representation, where $H_j$ is the
eigenspace of the character $\lambda\mapsto\lambda^j$. Then for all $i$ we
have  $F_i=\bigoplus_{j\leq i}H_j$, 
so that we have 
$$H=\gr H\,,$$
as graded algebras.  Here the grading on $H$ comes from the
$\cc^*$-representation, and the grading on $\gr H$ from the filtration on
$H$.\qed
\end{lem}

\begin{numex}\label{morestupex}
Consider the $\Gm$-action and equivariant vector field $W$ on $\pp^n$
of Example~\ref{stupex}. Then, because the zero locus $Z$ of $V$ has
dimension zero, the Koszul complex $K_W\com$ is a resolution of
$\O_Z$, and so we have that 
\begin{align*}
\hh^0(\pp^n,K_W\com)&=\Gamma(Z,\O_Z) \\
&=\cc[s_1,\ldots,s_n]/
(-s_1s_n,s_2-s_1^2,s_3-s_1s_2,\ldots,s_n-s_1s_{n-1})\\
&=\cc[s_1]/(-s_1^{n+1})\,.
\end{align*}
Note how the latter relations serve to recursively eliminate
$s_2,\ldots,s_n$, leaving only the first generator $s_1$ and the first
relation $-s_1s_n$.
\end{numex}

\begin{numrmk}\label{pqnz}
The assumption that $H^p(X,\Omega^q)=0$, for all $p\not=q$, was made
largely for simplicity, and because it is satisfied in all cases considered
in this paper.  It seems likely that the $E_1$-spectral sequence abutting
to $\hh(X,K_V\com)$ always degenerates. (See \cite{CaLie}, where this is
proved for the case of a compact K\"ahler manifold $X$ and a vector field
$V$ with non-empty zero set). If this is the case, then
$\hh^\ast(X,K_V\com)$ is a doubly graded $\cc$-algebra and isomorphic to
$H_{DR}(X)=\bigoplus_{p,q}H^p(X,\Omega^q)$ as such.
\end{numrmk}

\subsection{Equivariant actions on vector bundles}

To understand Chern classes  in the context of Proposition~\ref{AC}, we
need to study actions of vector fields on vector bundles.  First we recall
this concept without the presence of a $\Gm$-action.  Then we
consider the $\Gm$-equivariant case.

\subsubsection{Actions of vector fields on vector bundles}

Let $X$ be a smooth Deligne-Mumford stack and $E$ a vector bundle on $X$.

\begin{defn}[Carrell-Lieberman \cite{CL}]
Let $V$ be a vector field on $X$. An {\em action
}of $V$ on $E$ is a homomorphism of sheaves of $\cc$-vector spaces
$$\tilde{V}:E\longrightarrow E\,,$$
satisfying the Leibniz rule 
$$\tilde{V}(fe)=V(f)e+f\tilde{V}(e)\,,$$
for all local sections $f$ of $\O_X$ and $e$ of $E$. Here we interpret $V$
as a $\cc$-linear derivation $V:\O_X\to\O_X$.
\end{defn}

Note that any two actions of $V$ on $E$ differ by a homomorphism of vector
bundles $E\to E$. 

\begin{numrmk}\label{rein}
An action $\tilde{V}$ of $V$ on $E$ is the same thing as a $\Gm$-invariant
vector field $\ol{V}$ on $E$ lifting the vector field $V$ on $X$. Here we
mean $\Gm$-invariant with respect to the natural fiber-wise action of $\Gm$
on $E$ by scalar multiplication.  This means that in local
(linear) coordinates on $E$, the coefficients of $\ol{V}$ are linear
in these `vertical' coordinates.

The action $\tilde{V}$ is given in terms of the invariant lift
$\ol{V}$  by 
$$\tilde{V}(e)=De(V)-e\upst(\ol{V})\,,$$ 
for any local section $e:X\to E$. This is an equality of sections of
$e\upst\tT_E$. The bundle $e\upst\tT_E$ fits into the short exact sequence
$$\ses{E}{}{e\upst\tT_E}{}{\tT_X}\,,$$ 
which is canonically split by $De:\tT_X\to e\upst \tT_E$. We will often
identify the invariant lift $\ol{V}$ and the action $\tilde{V}$.

If the vector field $V$ on $X$ comes about by differentiating a
$\Ga$-action, and the vector field $\ol{V}$ on $E$ comes about by
differentiating a compatible {\em linear }$\Ga$-action on $E$, then
$\ol{V}$ is $\Gm$-invariant, and therefore gives rise to an action of $V$
on $E$.
\end{numrmk}

\begin{numexamples}\label{nablaex}
1\@. Let $\nabla:E\to\Omega_X\otimes E$ be a connection on $E$. Then, via
$\nabla$, every vector field $V$ acts on $E$. Just set $\tilde{V}$
equal to the covariant derivative $\nabla_V$.

2\@. Given the vector field $V$ on $X$, the vector bundles $\O_X$, $\tT_X$
and $\Omega_X$ have natural $V$-actions. For $\O_X$, we have $\tilde{V}=V$,
for $\tT_X$, we have that $\tilde{V}$ is equal to the Lie derivative with
respect to $V$ and for $\Omega_X$ we have $\tilde{V}=d\comp V$.

3\@. Given actions of $V$ on the vector bundles $E$ and $F$, there are
natural induced actions on $E\otimes F$ and $\sheafhom(E,F)$. These are
given by the usual Leibniz formulas
\begin{align*}
\tilde{V}(e\otimes f)& =\tilde{V}(e)\otimes f+e\otimes\tilde{V}(f)\\
\tilde{V}(\phi)(e)&
=\tilde{V}\big(\phi(e)\big)-\phi\big(\tilde{V}(e)\big)\,.
\end{align*}

4\@. Consider the vector field $W$ on $\pp^n$ from Example~\ref{stupex}.
It acts on $\O(1)$ by the formula
$$
\tilde{W}(x_i)=\begin{cases}x_{i+1}&\text{if $i<n$}\\
0&\text{if $i=n$}\,.\end{cases}
$$
We get induced actions of $W$ on $\O(m)$, for all $m\in\zz$. 
\end{numexamples}

\begin{numrmk}
Recall how connections on $E$ can be described in terms of splittings of
the short exact sequence of $\O_X$-modules (the Atiyah extension)
\begin{equation}\label{consplitsit}
\ses{\Omega_X\otimes E}{}{At(E)}{}{E}\,,
\end{equation}
where $At(E)$ is the $\O_X$-module whose underlying sheaf of $\cc$-vector
spaces is equal to $E\oplus(\Omega\otimes E)$ and whose $\O_X$-module
structure is defined by 
$$f\ast(e,\omega\otimes e')=(fe,f\omega\otimes e'+df\otimes e)\,.$$
If we denote by $s_0$ the $\cc$-linear splitting of
(\ref{consplitsit}) given by $e\mapsto (e,0)$, then every $\O_X$-linear
splitting $s$ of (\ref{consplitsit}) defines a connection on $E$ by the
formula $s=s_0+\nabla$, and conversely, every 
$\O_X$-linear splitting of (\ref{consplitsit}) comes from a unique
connection on $E$ in this way.

We can describe actions of $V$ on $E$ in a similar vein. We define a short
exact sequence of $\O_X$-modules
\begin{equation}\label{tildesplitsit}
\ses{E}{}{At_V(E)}{}{E}\,,
\end{equation}
where $At_V(E)$ is the $\O_X$-module whose underlying sheaf of $\cc$-vector
spaces is $E\oplus E$ and whose $\O_X$-modules structure is given by 
$$f\ast(e,e')=(fe,fe'+V(f)e)\,.$$
Moreover, the inclusion in (\ref{tildesplitsit}) is given by
$e'\mapsto(0,e')$ and the quotient map by $(e,e')\mapsto e$. Again, let us
denote the $\cc$-linear splitting $e\mapsto(e,0)$ of (\ref{tildesplitsit})
by $s_0$. Then actions of $V$ on $E$ and $\O_X$-linear splittings of
(\ref{tildesplitsit}) correspond bijectively to each other via the formula
$s=s_0+\tilde{V}$. 

Note that pushing out (\ref{consplitsit}) via the $\O_X$-homomorphism
$V\otimes\id_E:\Omega_X\otimes E\to E$ gives (\ref{tildesplitsit}). In
other words, we have a homomorphism of short exact sequences of
$\O_X$-modules
\begin{equation}\label{homses}
\vcenter{\xymatrix{
0\rto & \Omega_X\otimes E\dto_{V\otimes\id}\rto & At(E)\dto\rto &
E\dto^\id\rto 
& 0\phantom{\,.}\\ 
0\rto & E\rto & At_V(E)\rto & E\rto & 0\,.}}
\end{equation}
Thus every action $\tilde{V}$ of $V$ on $E$ gives rise to an $\O_X$-linear
map (the Carrell-Lieberman map)
\begin{align}\label{clvdef}
CL(\tilde{V}):At(E)&\longrightarrow E\\
(e,\omega\otimes e')&\longmapsto \langle\omega,V\rangle\nonumber
e'-\tilde{V}(e) \,,
\end{align}
making the diagram
\begin{equation*}
\vcenter{\xymatrix{
\Omega_X\otimes E\dto_{V\otimes\id}\rto & At(E) \dlto^{CL(\tilde{V})}\\
E&}}
\end{equation*}
commute.
\end{numrmk}

\subsubsection{The equivariant case}

Now suppose that $X$ is endowed with a $\Gm$-action, which has been lifted
to an action of $\Gm$ on $E$ by linear isomorphisms. Given a vector field
$V$ on $X$, and an action $\tilde{V}$ of $V$ on $E$, then
${^\lambda}\tilde{V}$ is an action of ${^\lambda}V$ on $E$.  So if $V$ is
equivariant, it is natural to consider equivariant lifts $\tilde{V}$, which
satisfy ${^\lambda}\tilde{V}=\lambda\tilde{V}$.

\begin{rmk}
Let $\ol{V}$ be an invariant vector field on $E$ (with respect to the
linear action), lifting the vector field $V$ on $X$. Let $\tilde{V}$ be the
associated action of $V$ on 
$E$. Then $\ol{V}$ is equivariant as a vector field on $E$ (with respect to
the geometric action)  if
and only if $\tilde{V}$ is equivariant as a $\cc$-linear homomorphism from
$E$ to $E$.

Hence, if the $\Gm$-action and the equivariant vector field $V$ on $X$ come
from an action of the group $G$ on $X$ as in Remark~\ref{pair}, then any
lift of the $G$-action to a vector bundle $E$ over $X$ gives rise to an
equivariant action of $V$ on $E$.
\end{rmk}

If $\tilde{V}$ and $V'$ are two equivariant actions of $V$ on $E$, then the
difference $h=\tilde{V}-V'$ is an equivariant homomorphism.

\begin{numex}
If $\nabla$ is an {\em invariant }connection on $E$, and $V$ an 
equivariant vector field on $X$, then the covariant derivative $\nabla_V$
is an 
equivariant action of $V$ on $E$.
\end{numex}

\begin{numex} \label{homotriv}
Suppose that $E$ is trivial and that there exists a trivialization of $E$
by a basis of global sections $(e_i)$, which are {\em homogeneous}, i.e., we
have 
$${^\lambda}e_i=\lambda^{r_i}e_i\,,$$
for all $i$ and certain integers $r_i$, with respect to the geometric
$\Gm$-action. Then the connection on $E$ induced by this trivialization via
the formula $\nabla(e_i)=0$ is invariant, and hence the covariant
derivative with respect to an equivariant vector field is equivariant.
\end{numex}

\begin{numex}\label{tildeW}
Consider the $\Gm$-action on $\pp^n$ from Example~\ref{stupex}.  It lifts
naturally to $\O(1)$ by the formula
$${^\lambda}x_i=\lambda^i x_i\,.$$
The action $\tilde{W}$ of $W$ on $\O(1)$ given by $\tilde{W}(x_i)=x_{i+1}$
($x_{n+1}=0$) is equivariant with respect to this $\Gm$-action. The same is
true for the induced actions on $\O(m)$, for all $m\in\zz$.
\end{numex}

\begin{numrmk}\label{justequiv}
Considering the natural induced action of $\Gm$ on the short exact
sequence~(\ref{consplitsit}), we note that all maps in (\ref{consplitsit})
are $\Gm$-invariant, and that invariant splittings correspond to invariant
connections.

Denote by $E^{(-1)}$ the vector bundle $E$ with the geometric $\Gm$-action
modified by the linear $\Gm$-action in such a way that
${^\lambda}(e^{(-1)})=\lambda {^\lambda}e$, where we have denoted by
$e^{(-1)}$ the section $e$ of $E$ considered as a section of $E^{(-1)}$. We
introduce the $\Gm$-action on~(\ref{tildesplitsit}) indicated by
$$\ses{E}{}{E^{(-1)}\oplus E}{}{E^{(-1)}}\,.$$
This choice of $\Gm$-action is necessary to ensure that the formula
$${^\lambda}\big(f\ast(e,e')\big)=({^\lambda}f)\ast{^\lambda}(e,e')$$
holds, for the $\O_X$-action on $At_V(E)$.

Now invariant splittings of~(\ref{tildesplitsit}) correspond to
equivariant actions of $V$ on $E$. Moreover, the homomorphism of short
exact sequences~(\ref{homses}) is of degree one (or equivariant, in our
language). The homomorphism $CL(\tilde{V}):At(E)\to E$ given by an
equivariant action of $V$ on $E$ is equivariant:
${^\lambda}\big(CL(\tilde{V})\big) = \lambda CL(\tilde{V})$.
\end{numrmk}

For future reference, we need some facts about the functorial behavior of
equivariant actions on vector bundles.

\begin{lem}[pullbacks]\label{tildepull}
Let $f:X\to Y$ be a morphism of smooth Deligne-Mumford stacks. Let $E$ be a
vector bundle over $Y$. Assume that $\Gm$ acts compatibly on $X$, $Y$ and
$E$. Let $V$ be an equivariant vector field on $X$ and $W$ an equivariant
vector field on $Y$, such that $Df(V)=W$. Let $\tilde{W}$ be an
equivariant action of $W$ on $E$. Then there is an induced
equivariant action $\tilde{V}$ of $V$ on $f\upst E$, such that 
$$\tilde{V}\big(f\upst(e)\big)=f\upst\big(\tilde{W}(e)\big)\,,$$
for every local section $e$ of $E$. If $V$, $W$ and $\tilde{W}$ come from
compatible actions of $G$ on $X$, $Y$ and $E$, then $\tilde{V}$ comes from
the induced $G$-action on $f\upst E$.
\end{lem}

\begin{lem}[pushforward]\label{tildepush}
Let $\pi:X\to Y$ be a flat and proper morphism of smooth Deligne-Mumford
stacks. Let $E$ be a vector bundle on $X$, such that $\pi\lst E$ is a
vector bundle on $Y$. Assume that $\Gm$ acts compatibly on $X$, $Y$ and
$E$. Let $V$ be an equivariant vector field on $X$ and $W$ an equivariant
vector field on $Y$, such that $D\pi(V)=W$. Let $\tilde{V}$ be an
equivariant action of $V$ on $E$. Then there is an induced equivariant
action $\tilde{W}$ of $W$ on $\pi\lst E$, such that 
$$\tilde{W}(e)=\tilde{V}(e)\,,$$
for every local section $e$ of $\pi\lst E$. If $V$, $W$ and $\tilde{V}$
come from compatible actions of $G$ on $X$, $Y$ and $E$, then $\tilde{W}$
comes from the induced $G$-action on $\pi\lst E$.
\end{lem}

\subsection{Chern classes}\label{Chern}

Now suppose that $X$ is a smooth Deligne-Mumford stack and $E$ a vector
bundle on $X$ of rank $r$. Suppose given compatible $\Gm$-actions on
$X$ and $E$. 
Finally, let $V$ be an equivariant vector field on $X$.

We shall now tensor the Koszul complex $K_V\com$ over
$\O_X$ with the sheaf of $\O_X$-modules $\sheafhom(E,E)$. Thus
$K_V\com\otimes\sheafhom(E,E)$ is a sheaf of differential graded modules
over the sheaf of differential graded $\O_X$-algebras $K_V\com$. Hence the
hypercohomology $\hh^0\big(X,K_V\com\otimes\sheafhom(E,E)\big)$ is a module
over the $\cc$-algebra $\hh^0(X,K_V\com)$. 

A $\cc^*$-representation on
$\hh^0\big(X,K_V\com\otimes\sheafhom(E,E)\big)$ is given by the composition
\begin{align*}
\hh^0\big(X,K_V\com\otimes\sheafhom(E,E)\big)
&\stackrel{\lambda^\ast}{\longrightarrow}
\hh^0\big(X,\lambda\upst K_V\com\otimes\lambda^\ast\sheafhom(E,E)\big)\\
&\stackrel{\psi_\lambda\otimes\rho}{\longrightarrow}
\hh^0\big(X,K_V\com\otimes\sheafhom(E,E)\big)\,,
\end{align*}
where $\psi_\lambda$ is the homomorphism given by~(\ref{psilam}) and
$\rho:\lambda\upst\sheafhom(E,E)\to\sheafhom(E,E)$ is the natural
isomorphism induced by the isomorphism $\lambda\upst E\cong E$ given by
the (geometric) action of $\Gm$ on $E$.

Note that the $\cc^*$ action on
$\hh^0\big(X,K_V\com\otimes\sheafhom(E,E)\big)$ is compatible with the
$\hh^0(X,K_V\com)$-action. Thus, via this $\cc^*$-representation,
$\hh^0\big(X,K_V\com\otimes\sheafhom(E,E)\big)$ becomes a graded
module over the graded $\cc$-algebra $\hh^0(X,K_V\com)$.

Now let $\tilde{V}$ be an equivariant action of $V$ on $E$. The
Carrell-Lieberman homomorphism $CL(\tilde{V})$ of~(\ref{clvdef}) gives
rise to a homomorphism of complexes
$$\xymatrix{
\ldots\rto& 0\dto\rto & \Omega_X\otimes E\dto^\id\rto& At(E)
\dto^{CL(\tilde{V})}\rto& 0\\ 
\ldots\rto & \Omega_X^2\otimes E\rto & \Omega\otimes E\rto^{V\otimes\id} &
E\rto &0}$$ 
which we can view as a homomorphism in the derived category from $E$ to
$K_V\com\otimes E$, because $[\Omega_X\otimes E\to At(E)]$ is a resolution
of $E$. We denote this homomorphism by 
$$c_{\tilde{V}}(E)\in\Hom_{D(\O_X)}(E,K_V\com\otimes E)\,.$$
Via the canonical identification
$$\Hom_{D(\O_X)}(E,K_V\com\otimes
E)=\hh^0\big(X,K_V\com\otimes\sheafhom(E,E)\big)\,,$$
$c_{\tilde{V}}(E)$ gives rise to a hypercohomology class, which we shall
also denote by 
$$c_{\tilde{V}}(E)\in\hh^0\big(X,K_V\com\otimes\sheafhom(E,E)\big)\,$$
and call the {\em Carrell-Lieberman class}.

It follows directly from Remark~\ref{justequiv} that the Carrell-Lieberman
class $c_{\tilde{V}}(E)$ is a degree one element of
$\hh^0\big(X,K_V\com\otimes\sheafhom(E,E)\big)$.

Now assume given a regular function $Q:M(r\times r)\to\aaa^1$ of degree
$p$, which is invariant under conjugation. This function gives rise to a
morphism of $X$-schemes
$$Q_E:\sheafhom(E,E)\longrightarrow\O_X\,,$$
by associating to an endomorphism $\phi$ of $E$ the number we get by
applying $Q$ to any matrix representation of $\phi$. The $X$-morphism $Q_E$
corresponds to a morphism of $\O_X$-algebras
$\O_X[t]\to\Sym\sheafhom(E,E)\dual$, where $t$ is a coordinate on
$\aaa^1$. Evaluating at $t$, this morphism gives rise to a global section
of $\Sym^p\sheafhom(E,E)\dual$, or equivalently, a symmetric $p$-linear
homomorphism 
$$Q':\Sym^p\sheafhom(E,E)\longrightarrow\O_X\,.$$
If $\phi$ is a local section of $\sheafhom(E,E)$, then the regular function
$Q_E\comp\phi$ on $X$ is equal to $Q'(\phi^p)$.

Recall how characteristic classes are defined in terms of the Atiyah
class. The Atiyah class $c(E)\in
H^1\big(X,\Omega\otimes\sheafhom(E,E)\big)$ is the cohomology class
corresponding to the extension of $\O_X$-modules $At(E)$ given
by~(\ref{consplitsit}) under the identification
$$\Ext^1(E,\Omega_X\otimes
E)=H^1\big(X,\Omega_X\otimes\sheafhom(E,E)\big)\,.$$
It gives rise to an element
$$c(E)^{\cup p}\in H^p\big(X,\Omega^{\otimes
p}\otimes\sheafhom(E,E)^{\otimes p}\big)\,,$$
by taking cup products.
Then we apply the map induced on $H^p$ by the composition
$$\Omega^{\otimes p}\otimes\sheafhom(E,E)^{\otimes p}
\longrightarrow\Omega^p\otimes \Sym^p\sheafhom(E,E)
\stackrel{\id\otimes Q'}{\longrightarrow}
\Omega^p\otimes\O_X=\Omega^p$$
to $c(E)^{\cup p}$. We obtain $c^Q(E)\in H^p(X,\Omega^p)$, the
characteristic class of $E$ defined by $Q$:
$$c^Q(E)=H^p(\id\otimes Q')\big(c(E)^{\cup p}\big)\,.$$

For example, if $Q$ is $(-1)^p$ times the degree $r-p$ coefficient of the
characteristic polynomial, then $c^Q(E)$ is the $p$-th Chern class of
$E$.

\begin{rmk}
If $X$ is a scheme, then under our identification of $H^p(X,\Omega^p)$ with
the singular cohomology $H^p(X,\cc)$, these Chern classes correspond to
the usual Chern classes.
\end{rmk}

Let us now apply a corresponding process to the Carrell-Lieberman class
$c_{\tilde{V}}(E)$. We start by taking the cup product of this class with
itself $p$ times: 
$$c_{\tilde{V}}(E)^{\cup p}\in\hh^0\big(X,(K_V\com)^{\otimes
p}\otimes\sheafhom(E,E)^{\otimes p}\big)\,.$$
Then we compose with the map induced on $\hh^0$ by 
$$(K_V\com)^{\otimes p}\otimes\sheafhom(E,E)^{\otimes p}
\stackrel{\mu\otimes Q'}{\longrightarrow}
K_V\com\otimes\O_X=K_V\com\,,$$
where $\mu:(K_V\com)^{\otimes p}\to K_V\com$ is the multiplication map
induced from the algebra structure on $K_V\com$. Thus we get the associated
characteristic class
\begin{equation}\label{cltoc}
c_{\tilde{V}}^Q(E)=\hh^0(\mu\otimes Q')\big(c_{\tilde{V}}(E)^{\cup p}\big)
\in\hh^0(X,K_V\com)\,.
\end{equation}

\begin{prop}
The class $c_{\tilde{V}}^Q(E)\in\hh^0(X,K_V\com)$ is homogeneous of degree
degree $p$.  If $X$ is proper and satisfies $h^{ij}(X)=0$, for all
$i\not=j$, then under the identification of Proposition~\ref{AC} we have
$$c_{\tilde{V}}^Q(E)=c^Q(E)\,.$$
\end{prop}
\begin{pf}
Let us denote by $K^{\geqslant-1}_V$ the na\"\i ve cutoff of $K_V\com$ at
minus 
one.  So $K^{\geqslant-1}_V$ denotes the two term complex
$[\Omega_X\to\O_X]$ given by the vector field $V$.  There is a 
canonical injection of complexes $K^{\geqslant-1}_V\to K_V\com$, which also
gives rise to an injection $K^{\geqslant-1}_V\otimes E\to K_V\com\otimes
E$. We also have a canonical projection $K^{\geqslant-1}_V\to\Omega_X[1]$,
giving rise to $K^{\geqslant-1}_V\otimes E\to\Omega_X\otimes E[1]$.  In the
derived category, we get an induced diagram
$$\xymatrix{
\Hom_{D(\O_X)}(E,K^{\geqslant-1}_V\otimes E)\dto\rto &
\Hom_{D(\O_X)}(E,K_V\com\otimes E)\\
\Hom_{D(\O_X)}(E,\Omega_X\otimes E[1]) &}$$
which we identify with the diagram
$$\xymatrix{
\Hom_{D(\O_X)}(E,K^{\geqslant-1}_V\otimes E)\dto\rto &
\hh^0\big(X,K_V\com\otimes\sheafhom(E,E)\big)\\
H^1\big(X,\Omega_X\otimes\sheafhom(E,E)\big) &}$$
Directly from the construction, it follows that the Carrell-Lieberman class
$c_{\tilde{V}}(E)$ lifts to $\Hom_{D(\O_X)}(E,K^{\geqslant-1}_V\otimes E)$,
and this lift maps to the Atiyah class:
$$\xymatrix{
CL(\tilde{V})\ar@{|->}[d]\ar@{|->}[r] &
c_{\tilde{V}}(E)\\
c(E) &}$$
The claim follows.
\end{pf}

\subsection{Localization to the big cell}

Let $X$ be, as above, a
smooth Deligne-Mumford stack with a $\Gm$-action, and let $E$ be a vector
bundle to which the $\Gm$-action has been lifted. Let $V$ be an equivariant
vector field on $X$ and $\tilde{V}$ an equivariant action of $V$ on
$E$. Finally, let us denote by $Z\subset X$ the closed substack defined by
the vanishing of $V$: the structure sheaf $\O_Z$ of $Z$ is defined to be
the cokernel of the vector field $V:\Omega_X\to\O_X$. Another way to think
of $\O_Z$ is as the zero-degree cohomology sheaf of $K_V\com$. 

There is a canonical morphism of sheaves of differential graded algebras
$K_V\com\to\O_Z$, inducing a canonical morphism of graded algebras
$$\hh^0(X,K_V\com)\longrightarrow \Gamma(Z,\O_Z)\,.$$

\begin{rmk}
When we restrict to $Z$, the action $\tilde{V}$ induces an $\O_Z$-linear
homomorphism $\tilde{V}_Z:E_Z\to E_Z$. Applying the invariant polynomial
$Q$, we get a regular function 
$$Q(\tilde{V}_Z)\in\Gamma(Z,\O_Z)\,.$$
Under the canonical morphism
$\hh^0(X,K_V\com)\rightarrow \Gamma(Z,\O_Z)$ the characteristic class
$c_{\tilde{V}}^Q(E)$ maps to $Q(\tilde{V}_Z)$. 
\end{rmk}

If $X$ is affine and $E$ is trivial over $X$, trivialized by a
homogeneous basis $(e_i)$ as in Example~\ref{homotriv}, then we form the
matrix $M$ of $\tilde{V}$ with respect to this basis. In other words,
$M=(m_{ij})$ is a square matrix with entries $m_{ij}\in\Gamma(X,\O_X)$,
characterized by
$$\tilde{V}(e_i)=\sum_j m_{ji} e_j\,.$$
In this case we have 
$$\Gamma(Z,\O_Z)=\Gamma(X,\O_X)/V\big(\Gamma(X,\Omega_X)\big)$$
and we can compute $Q(\tilde{V}_Z)\in\Gamma(Z,\O_Z)$ simply as 
the congruence class of $Q(M)\in\Gamma(X,\O_X)$.

Note that if $e_i$ is homogeneous of degree $d_i$, then $m_{ij}$ is
homogeneous of degree $d_i-d_j+1$.

\begin{example}
Returning to Example~\ref{morestupex}, using the equivariant action
$\tilde{W}$ of $W$ on $\O(m)$ of Example~\ref{tildeW}, we see that a
basis for $\O(m)$ over $\{x_0=1\}\subset\pp^n$ is given by $x_0^m$. We
have $\tilde{W}(x_0^m)=m\frac{x_1}{x_0} x_0^m=m s_1 x_0^m$ in the
notation of Example~\ref{morestupex}. Thus the matrix of $\tilde{W}$
with respect to this basis is $m s_1$ and so the first Chern class
$c_1\big(\O(m)\big)\in\hh^0(\pp^n,K_W\com)$ is equal to $ms_1$.
\end{example}

We shall apply these ideas in the following context. The stack $X$ will be
proper, smooth and satisfy our assumption on the vanishing of off-diagonal
Hodge numbers.  The stack $X$ will be endowed with a $\Gm$-action and
several equivariant vector bundles $E_i$. Moreover, $X$ will have an
equivariant vector field $V$ on it, with natural equivariant lifts
$\tilde{V}_i$ to the various $E_i$.  We will construct an affine scheme
$T$, with a lift of the $\Gm$-action to $T$ and an \'etale morphism $T\to
X$. Over $T$ we trivialize all $E_i$ 
using homogeneous bases.  Finally, we choose a collection of invariant
polynomials $Q_i$ of various degrees, giving rise to characteristic classes
$c_i=c_i^{Q_i}(E_i)\in H_{DR}(X)$.

We consider the natural morphism
\begin{equation}\label{c9}
H_{DR}(X)=\hh^0(X,K_V\com)  \longrightarrow
\hh^0(T,K_V\com) =
\Gamma(T,\O_T)/V\big(\Gamma(T,\Omega_T)\big)\,,
\end{equation}
and note that we can compute the images of $c_i$ under~(\ref{c9}) as
$Q_i(M_i)$, where $M_i$ is the matrix representation of $\tilde{V}_i$ with
respect to our homogeneous basis of $E_i\resto T$.

In the cases we consider, it will turn out that the map~(\ref{c9}) is
injective in sufficiently low degrees and that $H_{DR}(X)$ is
generated (as a $\cc$-algebra) by the classes $c_i$.


\section{Preliminaries on stable maps}\label{prelims}

\subsection{Stable maps to vector bundles}

Let $X$ be a smooth and proper algebraic variety and $E$ a vector bundle
over $X$, with structure morphism $p:E\to X$.

\begin{note}
Pulling back via $p$ induces an isomorphism of Picard groups
$p\upst:\Pic(X)\to \Pic(E)$ which preserves the ample line bundles. Hence
we have a canonical isomorphism $H_2(E)^+\to H_2(X)^+$ which we use to
identify these two semi-groups.
\end{note}

\begin{note}
Let $\tilde{f}:C\to E$ be a morphism from a prestable marked curve $(C,x)$
to $E$. Then $\tilde{f}$ is stable of class $\beta$ if and only if
$p(\tilde{f})=p\comp \tilde{f}$ is stable of class $\beta$. This is because
$p( \tilde{f})$ cannot contract any component of $C$ which $\tilde{f}$
does not already contract, because no component of $C$ can map into a fiber
of $p$ , these fibers being affine.
\end{note}

Let 
$$\xymatrix{\cC\dto_{\pi}\rto^{f} & X\\ 
\ol{M}_{g,n}(X,\beta) &}$$
be the universal curve and universal stable  map.  We may view $\pi\lst
f\upst E$ as a (relative) scheme over $\ol{M}_{g,n}(X,\beta)$. Then it
represents the following functor:
$$\pi\lst f\upst E(T)=\Gamma(\cC_T,f\upst E)\,,$$
for any $\ol{M}_{g,n}(X,\beta)$-scheme $T$.  Here $\cC_T$ abbreviates the
pull-back of $\cC$ to $T$. 
If $T\to\ol{M}_{g,n}(X,\beta)$ is given by $(C,x,f)$, then
\begin{align*}
\pi\lst f\upst E(T)&=\Gamma(C, f\upst E)\\
&=\{\tilde{f}:C\to E\st p(\tilde{f})=f\}\,.
\end{align*}
Hence $\pi\lst f\upst E$ as a stack in its own right (forgetting the
$\ol{M}_{g,n}(X\beta)$-structure) associates to the $\cc$-scheme $T$ the
groupoid of triples $(C,x,\tilde{f})$, where $(C,x)$ is a prestable marked
curve and $\tilde{f}:C\to E$ is a morphism such that $p(\tilde{f}):C\to X$
is a stable map of class $\beta$.  by the above notes, this is equivalent
to saying that $\tilde{f}$ is a stable map of class $\beta$. Thus we
conclude that 
$$\pi\lst f\upst E=\ol{M}_{g,n}(E,\beta)\,.$$
In particular, $\ol{M}_{g,n}(E,\beta)$ is an algebraic stack, representable
over $\ol{M}_{g,n}(X,\beta)$, by an (abelian) cone (see \cite{BF} for
this terminology). We note the following
consequence:

\begin{prop}\label{VB}
If $E$ is convex over $X$, i.e., $H^1(\pp^1,f\upst E)=0$ for all morphisms
$f:\pp^1\to X$, then $\ol{M}_{0,n}(E,\beta)$ is a vector bundle over
$\ol{M}_{0,n}(X,\beta)$, canonically identified $\pi\lst f\upst E$. The
rank of this vector bundle is $\langle c_1(E),\beta\rangle+\rk E$.
\end{prop}

\begin{numex}
For $m\geq0$, the vector bundle $E_m=\pi\lst f\upst\O_{\pp^n}(m)$ over 
$\ol{M}_{0,\nu}(\pp^n,d)$ has rank $md+1$ and represents
$\ol{M}_{0,\nu}(\O_{\pp^n}(m),d)$. 
\end{numex}

\subsection{Stable maps to $\pp^\infty$}
\label{sminf}

Let $n<m$ be integers and consider $\pp^n$ as a subvariety of $\pp^m$
via $\pp^n=\{\langle x_0,\ldots,x_n,0,\ldots,0\rangle\}\subset\pp^m$. 
Let $H=\langle
0,\ldots,0,x_{n+1},\ldots,x_m\rangle\}\cong\pp^{m-n-1}$.  Let $U=\pp^m-H$
and consider the 
projection with center $H$ onto $\pp^n$, which is defined on $U$ and makes
$U$ a vector bundle over $\pp^n$ of rank $m-n$, in fact this vector bundle
is isomorphic to a direct sum of $m-n$ copies of
$\O_{\pp^n}(1)$.

By applying $\ol{M}_{0,\nu}(\argument,d)$ to the diagram of varieties
$$\xymatrix{
U\dto\ar@{^{(}->}[r] & \pp^m \\
\pp^n & }$$
we get the diagram of stacks 
$$\xymatrix{
\ol{M}_{0,\nu}(U,d)\dto_{\rho}\ar@{^{(}->}[r]^{\iota} &
\ol{M}_{0,\nu}(\pp^m,d) \\ 
\ol{M}_{0,\nu}(\pp^n,d) &.}$$
Here $\rho$ is a vector bundle of rank $(d+1)(m-n)$. 
Let $\kappa$ be the zero section of this vector bundle. The projection
$\rho$ is a homotopy equivalence implying that
$$\kappa\upst:H^p\big(\ol{M}_{0,\nu}(U,d)_\an,\cc\big)\longiso
H^p\big(\ol{M}_{0,\nu}(\pp^n,d)_\an,\cc\big)$$
is an isomorphism for all $p$. 

The complement of $\ol{M}_{0,\nu}(U,d)$ in $\ol{M}_{0,\nu}(\pp^m,d)$
consists of all stable maps to $\pp^m$ whose image intersects $H$. The
locus of these stable maps has codimension $n$ in
$\ol{M}_{0,\nu}(\pp^m,d)$. Thus
$$\iota\upst:H^p\big(\ol{M}_{0,\nu}(\pp^m,d)_\an,\cc\big)
\longrightarrow H^p\big(\ol{M}_{0,\nu}(U,d)_\an,\cc\big)$$
is an isomorphism for all $p\leq 2n-2$, by cohomological purity.

We conclude that
\begin{equation}\label{stabilizeagain}
(\iota\kappa)\upst:
H^p\big(\ol{M}_{0,\nu}(\pp^m,d)_\an,\cc\big)\longrightarrow 
H^p\big(\ol{M}_{0,\nu}(\pp^n,d)_\an,\cc\big)
\end{equation}
is an isomorphism for $p\leq2n-2$. 
This leads to the following definition: 

\begin{defn}
For every $p\geq0$ we define
$$H^p\big(\ol{M}_{0,\nu}(\pp^\infty,d)_\an,\cc\big)=
H^p\big(\ol{M}_{0,\nu}(\pp^n,d)_\an,\cc\big)\,,$$
for any $n$ such that $n\geq\frac{1}{2}(p+2)$.
\end{defn}

It follows from the above considerations that any two different choices of
$n$ lead to canonically isomorphic definitions of
$H^p\big(\ol{M}_{0,\nu}(\pp^\infty,d)_\an,\cc\big)$, for fixed $p$.
Taking the direct sum over all $p$ gives rise to the $\cc$-algebra
$$H\upst\big(\ol{M}_{0,\nu}(\pp^\infty,d)_\an,\cc\big)\,.$$

For any $n$ we have the canonical restriction map
$$H\upst\big(\ol{M}_{0,\nu}(\pp^\infty,d)_\an,\cc\big)\longrightarrow
H\upst\big(\ol{M}_{0,\nu}(\pp^n,d)_\an,\cc\big)$$
which is a $\cc$-algebra morphism and an isomorphism in degrees less
than $n$ (recall our degree convention from Remark~\ref{allthecoh}). 

\begin{numrmk}
For every $\nu$, $d$, $n$ we have
$$H\upst\big(\ol{M}_{0,\nu}(\pp^n,d)_\an,\cc\big)=
H_{DR}\big(\ol{M}_{0,\nu}(\pp^n,d)\big)=
\bigoplus_{p}H^p\big(\ol{M}_{0,\nu}(\pp^n,d),\Omega^p\big)\,.$$
The first equality follows from Remark~\ref{allthecoh}. The second equality
is proved using the technique of virtual Hodge
polynomials to reduce to the strata of a suitable stratification. Use
the stratification by topological type. The details are worked out by
Getzler and Pandharipande in~\cite{getzler}. This justifies writing also
$$H\upst\big(\ol{M}_{0,\nu}(\pp^\infty,d)_\an,\cc\big)=
H_{DR}\big(\ol{M}_{0,\nu}(\pp^\infty,d)\big)=
\bigoplus_{p}H^p\big(\ol{M}_{0,\nu}(\pp^\infty,d),\Omega^p\big)\,.$$
\end{numrmk}

\begin{numrmk} \label{remarksurjective}
It is, in fact, true that
$$
(\iota\kappa)\upst:
H^p\big(\ol{M}_{0,\nu}(\pp^m,d)_\an,\cc\big)\longrightarrow 
H^p\big(\ol{M}_{0,\nu}(\pp^n,d)_\an,\cc\big)
$$
is {\em surjective }for all $p$. One way to prove this is as
follows. Consider the $\Gm$-action on $\ol{M}_{0,\nu}(\pp^m,d)$
induced by the $\Gm$-action on $\pp^m$ given by
$$\langle x_0,\ldots,x_n,x_{n+1},\ldots,x_m\rangle\cdot\lambda
=\langle x_0,\ldots,x_n,\lambda x_{n+1},\ldots,\lambda x_m\rangle\,.$$
The substack $\ol{M}_{0,\nu}(\pp^n,d)$ is a fixed locus for this
action and $\ol{M}_{0,\nu}(U,d)$ is the big  Bia\l ynicki-Birula cell
of this action. Using virtual Poincar\'e polynomials as in
\cite{getzler} one proves that all  Bia\l ynicki-Birula cells have
only even cohomology, which implies the claim. (All fixed loci can be
described explicitly, showing that they are amenable to the techniques
of \cite{getzler}. See for example \cite{kontsevich} or
\cite{cetraro} for explicit descriptions of fixed loci of
$\Gm$-actions on stable map stacks.) 

This result implies that
$$H\upst\big(\ol{M}_{0,\nu}(\pp^\infty,d)_\an,\cc\big)\longrightarrow
H\upst\big(\ol{M}_{0,\nu}(\pp^n,d)_\an,\cc\big)$$
is an epimorphism of $\cc$-algebras for all $n$.

As we do not have a reference for the Bia\l ynicki-Birula
decomposition of a Deligne-Mumford stack, we are careful to point out
where we use this result. The only places are, in fact,
Corollary~\ref{generation}, Proposition~\ref{evidence},
Corollary~\ref{td1} and Proposition~\ref{td2}.
\end{numrmk}

\section{Parameterizing stable maps to $\pp^n$}\label{SecPar}

Let $d$ be an integer greater than or equal to $3$.
We will prove that there is an open substack $U$ of $\ol{M}_{0,0}(\pp^n,d)$
which is a vector bundle over $[\ol{M}_{0,d}/S_d]$. In fact, $U=[T/S_d]$,
where $T$ is a vector bundle over the scheme $\ol{M}_{0,d}$. The whole
stack $\ol{M}_{0,0}(\pp^n,d)$ can be covered by substacks isomorphic to
$U$. 

One way to describe $U$ is as follows: introduce on $\pp^n$ a suitable
$\Gm$-action and consider the induced $\Gm$-action on
$\ol{M}_{0,0}(\pp^n,d)$. Then consider to every fixed component $U_0$ of this
$\Gm$-action the associated substack $U$ of all points that move to $U_0$
as $\lambda\to0$, $\lambda\in\Gm$. For a suitable component $U_0$, we have
that $U$ is open in $\ol{M}_{0,0}(\pp^n,d)$ and that $U_0$ is isomorphic to
$[\ol{M}_{0,d}/S_d]$. 

The most important case is the case where $n=1$. In this case $U$ can
simply be described as the stack of all stable maps unramified over
$\infty\in\pp^1$. The case of general $n$ is easily reduced to this special
case.

We start by defining a vector bundle $T$ over $\ol{M}_{0,d}$ and
constructing a stable map $(\tilde{C},f)$ to $\pp^1$, parametrized by $T$.

\subsection{For every $i$ a degree 1 map to $\pp^1$}
\label{forevery}

Let $T$ be a scheme and $(C,x)$ a prestable curve of genus zero over $T$,
with $d$ marked points.
We denote the structure map by $\pi:C\to T$ and $x$ stands for the
$d$-tuple of sections $x_i:T\to C$, $i=1,\ldots,d$. Consider the
canonical line bundles
$$\omega_i=x_i\upst\Omega_{C/T}$$ on $T$ and their duals
$L_i=\omega_i\dual$. Let $D_i\hookrightarrow C$ be the effective Cartier
divisor defined by the $i$-th section $x_i$.  By $\kK_C$ we denote the
sheaf of total quotient rings of $C$ and by $\kK\upst_C$ its sheaf of units.

Note that the sheaf $\pi\lst\O(D_i)$ is locally free of rank 2,
because the genus of $C/T$ is 0.

We will now fix an index $i$ and define a canonical homomorphism
$$h_i:L_i\longrightarrow\pi\lst\O(D_i).$$
Note that $h_i$ may equivalently be defined by a global section
$$h_i\in\Gamma(C,\pi\upst\omega_i(D_i)).$$

Assume for the moment that there exists a global section
$s\in\Gamma(C,\kK_C\upst)$, which in a neighborhood of $D_i$ generates
the ideal sheaf of $D_i$, and which does not vanish anywhere else,
except at $D_i$.
In other words, $s$ is the reciprocal of a global section of
$\pi\lst\O(D_i)$ which is nowhere contained in the submodule
$\pi\lst\O=\O$.
Then we define $h_i$ by the formula
\begin{equation}\label{formhi}
h_i:=\left(\frac{1}{s}+\frac{1}{1-d}\sum_{j\not=i}\frac{1}{s(x_j)}\right)
      ds(x_i)\,\,.
\end{equation}
Note that this makes sense, even if $s$ has a pole at any of the
$x_j$.

\begin{lem}
Equation~(\ref{formhi}) defines $h_i$ independently of the choice of
$s$.
\end{lem}
\begin{pf}
Let $t$ be another section of $\kK_C\upst$, satisfying the same
conditions as $s$. Then $\frac{t}{s}\in\kK\upst$ is an element of
$\O\upst$ in a neighborhood of $D_i$ and so we have
\begin{align*}
dt(x_i)&=(d\textstyle\frac{t}{s}s)(x_i)\\
       &=s(x_i)d\textstyle\frac{t}{s}(x_i)+\frac{t}{s}(x_i)ds(x_i)\\
       &=\frac{t}{s}(x_i)ds(x_i)\,.
\end{align*}
Note also that
$$\frac{1}{s}-\frac{1}{t}\phantom{\cdot}\frac{t}{s}(x_i)$$
is a {\em regular }function on $C$, since the only poles cancel out.
Thus this regular function is constant on the fibers of $\pi$ and so
it is equal to its evaluation at any of the sections $x_j$. In other
words,
$$\frac{1}{s}-\frac{1}{t}\phantom{\cdot}\frac{t}{s}(x_i)=
\frac{1}{s(x_j)}-\frac{1}{t(x_j)}\phantom{\cdot}\frac{t}{s}(x_i)\,\,,$$
which, again, also makes sense if $s$ or $t$ has a pole at $x_j$.

We may now calculate as follows:
\begin{align*}
\begin{split}
&\left(\frac{1}{t}+\frac{1}{1-d}\sum_{j\not=i}\frac{1}{t(x_j)}\right)dt(x_i)\\
&\hskip1cm=\left(\frac{1}{t}+\frac{1}{1-d}\sum_{j\not=i}\frac{1}{t(x_j)}\right)
  \frac{t}{s}(x_i)ds(x_i)
\end{split}\\
&\hskip1cm=\left(\frac{1}{s}-\left(\frac{1}{s}-\frac{1}{t}\phantom{\cdot}
  \frac{t}{s}(x_i)\right)
  +\frac{1}{1-d}\sum_{j\not=i}\frac{1}{t(x_j)}\phantom{\cdot}
  \frac{t}{s}(x_i)\right)   ds(x_i) \\
&\hskip1cm=\left(\frac{1}{s}
  +\frac{1}{1-d}\sum_{j\not=i}\left(\frac{1}{s}-\frac{1}{t}\phantom{\cdot}
  \frac{t}{s}(x_i)+\frac{1}{t(x_j)}\phantom{\cdot}
  \frac{t}{s}(x_i)\right)\right)   ds(x_i) \\
&\hskip1cm=\left(\frac{1}{s}
  +\frac{1}{1-d}\sum_{j\not=i}
  \left(\frac{1}{s(x_j)}-\frac{1}{t(x_j)}\phantom{\cdot}
  \frac{t}{s}(x_i)+\frac{1}{t(x_j)}\phantom{\cdot}
  \frac{t}{s}(x_i)\right)\right)   ds(x_i)\\
&\hskip1cm=\left(\frac{1}{s}
  +\frac{1}{1-d}\sum_{j\not=i}
  \frac{1}{s(x_j)}\right)   ds(x_i)\,\,.
\end{align*}
Thus, indeed, $h_i$ is well-defined.
\end{pf}

\begin{cor}\label{hi}
There exists a unique homomorphism
$$h_i:L_i\longrightarrow\pi\lst\O(D_i)\,\,,$$
such that the restriction of $h_i$ to any open subset of $T$ which
admits an $s$ as above is given by Formula~(\ref{formhi}). The inverse
image under $h_i$ of the submodule $\pi\lst\O\subset\pi\lst\O(D_i)$ is
$0$.
\end{cor}
\begin{pf}
Since Zariski-locally in $T$ we can find an $s$ as
required, these locally defined $h_i$ glue.
\end{pf}

The basic properties of $h_i$ are summarized in the following

\begin{prop}\label{firstchar}
Evaluating $h_i$ at $x_j$, for $j\not=i$ defines canonical sections
$h_i(x_j)\in\Gamma(T,\omega_i)$. We have
$$\sum_{j\not=i}h_i(x_j)=0\,.$$

Evaluating $h_i$ at $x_i$ gives a canonical section of
$\Gamma(T,\omega_i\otimes x_i\upst\O(D_i)=\Gamma(T,\O_T)$. We have
$$h_i(x_i)=1\,,$$ under this identification.

Finally, $h_i$ is characterized completely by these two properties.
\end{prop}
\begin{pf}
The two properties mentioned follow directly from the explicit definition
of $h_i$ in terms of a local parameter given above. The fact that $h_i$ is
determined by two properties follows from the fact that $\pi\lst\O(D_i)$ is
of rank two.
\end{pf}

Let us now suppose given a section $\tau\in\Gamma(T,L_i)$. Then
$h_i(\tau)\in\Gamma(T,\pi\lst\O(D_i))=\Gamma(C,\O(D_i))$ is a
meromorphic function on $C$, which we may view as a rational map
$$\xymatrix{h_i(\tau):C\ar@{.>}[r]&\pp^1}\,\,.$$
This rational map is defined by the pencil given by the global
sections $1,h_i(\tau)$ of $\O(D_i)$ (at least in the case where these
global sections define a pencil). Note that if $\tau$ is nowhere
vanishing, then $h_i(\tau):C\to\pp^1$ is an everywhere defined
morphism. 

For simplicity, assume now that $T$ is smooth and that $V=\{\tau=0\}$
is a Cartier divisor. Then $V_C=\pi^{-1}V$ is also a Cartier
divisor. Consider the closed subscheme $Z_i=D_i\cap V_C$ of
$C$. \'Etale locally in $C$, we can find coordinates
$(t_1,\ldots,t_r)$ for $T$ and $(t_1,\ldots,t_r,s)$ for $C$ such that
$D_i$ is given by $s=0$ and $V_C$ by $t_1=0$.

Thus the structure of the blow-up $\tilde{C}$ of $C$ along
$Z_i=\{s=t_1=0\}$ is transparent: $\tilde{\pi}:\tilde{C}\to T$ is again
a family of genus 0 prestable curves. The divisors $D_j$, for
$j\not=i$ are contained in the locus where $\tilde{C}\to C$ is an
isomorphism.  Thus they are divisors on $\tilde{C}$ and are images of
sections $x_j:T\to\tilde{C}$. The strict transform $\tilde{D}_i$ of
$D_i$ is the image of another section $\tilde{x}_i:T\to\tilde{C}$. The
pair $(\tilde{C},\tilde{x})$, where $\tilde{x}_j=x_j$, for $j\not=i$,
is a
prestable marked curve. Moreover, the blow up morphism $p:\tilde{C}\to
C$ is a morphism of prestable curves, in particular,
$p(\tilde{x}_i)=x_i$. 

Let $\tilde{L}_i=\tilde{\omega}_i\dual$, where
$\tilde{\omega}_i=\tilde{x}_i\upst\Omega_{\tilde{C}/T}$. There is a
natural exact sequence
$$\ses{\tilde{L}_i}{}{L_i}{}{\O_V}\,\,,$$
coming  via ${\tilde{x}_i}\upst$ from the exact sequence 
$$\ses{\omega_{\tilde{C}}\dual}{}{p\upst\omega_C}{}{\O_E}\,,$$
where $E\subset \tilde{C}$ is the exceptional divisor. Thus
$\tilde{L}_i=L_i(-V)$.  The section $\tau$ of $L_i$ factors
through the subsheaf $\tilde{L}_i$; let us denote this section of
$\tilde{L}_i$ by $\tilde{\tau}$. The section $\tilde{\tau}$ is nowhere
vanishing, so it trivializes $\tilde{L}_i$. 

The marked prestable curve $(\tilde{C},\tilde{x})$ has the canonical
morphism 
$$\tilde{h}_i:\tilde{L}_i\longrightarrow\tilde{\pi}\lst\O(\tilde{D}_i)$$
associated to it.  Since $\tilde{\tau}$ is nowhere vanishing, the
associated meromorphic function $\tilde{h}_i(\tilde{\tau})$ defines a
morphism
$$\tilde{h}_i(\tilde{\tau}):\tilde{C}\longrightarrow\pp^1.$$

The base locus of the pencil defined by the global sections
$1,h_i(\tau)$ of $\O_C(D_i)$ is equal to $Z_i=D_i\cap V_C$. The
rational map $\xymatrix{h_i(\tau):C\ar@{.>}[r]&\pp^1}$ sends the
divisor $D_i$ to $\infty\in\pp^1$ and the divisor $V_C$ to
$0\in\pp^1$.

We have a commutative diagram
$$\xymatrix{
\tilde{L}_i\dto_{\tilde{h}_i}\rrto && L_i\dto^{h_i}\\
\tilde{\pi}\lst\O(\tilde{D}_i)\drto && \pi\lst\O(D_i)\dlto\\
&\pi\lst\kK_C&}$$
Thus the morphism
$\tilde{h}_i(\tilde{\tau}):\tilde{C}\rightarrow\pp^1$ is the morphism
defined by the rational map $\xymatrix{h_i(\tau):C\ar@{.>}[r]&\pp^1}$
via blowing up the locus where it is undefined.

To abbreviate notation, let us write $f$ for the morphism
$\tilde{h}_i(\tilde{\tau}):\tilde{C}\rightarrow\pp^1$.

Before we can state the next proposition about $f$, we need to recall
a few facts about ramification.  Let $\pi:C\to S$ be a family of
prestable curves over a scheme $S$ and let $f:C\to\pp^1$ be a
morphism. Then there is a canonical homomorphism of line bundles on $C$
$$f\upst\Omega_{\pp^1}\longrightarrow\omega_{C/S}\,\,,$$
where $\omega_{C/S}$ is the relative dualizing sheaf of $C$ over $S$
(which is also equal to the determinant of the relative cotangent
complex, which shows the existence of the
homomorphism).
The ideal sheaf
$$I=( f\upst\Omega_{\pp^1}: \omega_{C/S})= \{a\in\O_C\st
a\cdot\omega_{C/S}\subset f\upst\Omega_{\pp^1}\}$$
defines a closed subscheme $R\subset C$ called the {\em ramification
scheme }of $f$. There is a canonical exact sequence
\begin{equation}\label{notses}
f\upst\Omega_{\pp^1} \longrightarrow\omega_{C/S}\longrightarrow
\omega_{C/S}\otimes\O_R\longrightarrow0\,\,.
\end{equation}
We say that $f$ is {\em unramified }over $a\in\pp^1$, if
$f^{-1}(a)\cap R=\varnothing$. Note that this is an open condition in
$S$. Moreover, $f$ is unramified over $a\in\pp^1$, if and only if
$(f\upst\Omega_{\pp^1}\to\omega_{C/S})\resto f^{-1}(a)$ is
surjective. This means that for every geometric point $s$ of $S$, the
map $f_s:C_s\to\pp^1$ is unramified over $a$. In particular,
$f_s^{-1}(a)$ consists of $\deg f_s$ distinct points.

\begin{lem}\label{unramlem}
Let $f$ be unramified over $a\in\pp^1$. Then $f^{-1}(a)$ is finite
\'etale of degree $\deg f$ over $S$.  If $x:S\to C$ is a section such
that $f(x)=a$, then the derivative 
$$Df:\tT_{C/S}\longrightarrow f\upst\tT_{\pp^1}$$
pulls back via $x$ to an isomorphism
$$Df(x):x\upst\tT_{C/S}\longrightarrow
\tT_{\pp^1}(a)\tensor\O_S\,\,,$$ 
where $\tT_{\pp^1}(a)$ denotes the tangent space of $\pp^1$ at 
$a$. 
\end{lem}
\begin{pf}
The structure sheaf $\O_{f^{-1}(a)}$ of the inverse image
$f^{-1}(a)$ fits into the exact sequence 
$$\ses{\O_C}{}{f\upst\O(1)}{}{\O_{f^{-1}(a)}}\,,$$
obtained from an identification $\O(1)=\O(a)$. 
We get an induced exact sequence
$$\ses{\O_S}{}{\pi\lst f\upst \O(1)}{}{\pi\lst\O_{f^{-1}(a)}}\,.$$
Since this sequence stays exact after arbitrary base change, 
it proves that $\pi\lst\O_{f^{-1}(a)}$ is locally free of rank
$\deg f$ (recall that $\pi\lst f\upst\O(1)$ is locally free of rank
$\deg f+1$). By considering the fibers of $C\to S$ we see that
$f^{-1}(a)$ is quasi-finite over $S$. Since it is also proper, it is
finite, thus flat. Then \'etale is equivalent to unramified, which can
be checked fiber-wise.
\end{pf}

Now we come back to $f:\tilde{C}\to\pp^1$ over $T$, defined above.

\begin{prop}\label{propdeg1}
The morphism $f:\tilde{C}\to\pp^1$ is a family of degree 1 maps to
$\pp^1$, unramified over $\infty\in\pp^1$. 
We have
$$f^{-1}(\infty)=\tilde{D}_i\,\,,$$
and 
\begin{equation}\label{ford}
Df(\tilde{x}_i)(\tilde{\tau})= \frac{\partial}{\partial
z}(\infty)\otimes 1\,\,,
\end{equation}
where $z$ is the canonical coordinate at $\infty\in\pp^1$ and
$\frac{\partial}{\partial z}(\infty)$ is the evaluation of
$\frac{\partial}{\partial z}$ at $z=0$.
\end{prop}
\begin{pf}
This follows directly from the  construction. Formula~(\ref{ford}) is
a direct calculation.
\end{pf}

\subsection{A degree $d$ map to $\pp^1$}\label{degn}

Now let us suppose given a regular function $b\in\Gamma(T,\O)$ and for
every $i=1,\ldots,d$ a section
$\tau_i\in\Gamma(T,L_i)$. Then for each $i$ we have
$h_i(\tau_i)\in\Gamma(C,\O(D_i))$ and so for the sum we have
$$b+\sum_{i=1}^d
h_i(\tau_i)\in\Gamma\big(C,\O(\textstyle\sum_{i}D_i)\big)\,\,,$$
which we may also view as a rational map
\begin{equation}\label{ratmap}
\xymatrix{b+\sum_{i=1}^d h_i(\tau_i):C\ar@{.>}[r]&\pp^1}\,\,.
\end{equation}
Note that this is an everywhere defined morphism, if all of the
$\tau_i$ are nowhere vanishing.

Let $V_i=\{\tau_i=0\}$ and, as above,
$$Z_i=D_i\cap V_{i,C}\,\,.$$
Let also $Z=Z_1\cup\ldots\cup Z_d$ be the union of these pairwise
disjoint closed subschemes.

Let us assume, as above, that $T$ is smooth and that the $V_i$ are
Cartier divisors. Let $\tilde{C}$ be the blow-up of $C$ along $Z$. We
get induced sections $\tilde{x}_i:T\to\tilde{C}$, making
$(\tilde{C},\tilde{x})$ a prestable marked curve. We also get induced
nowhere vanishing sections $\tilde{\tau}_i\in\Gamma(T,\tilde{L}_i)$ 
and hence an everywhere defined morphism
$$f=b+\sum_{i=1}^d \tilde{h}_i(\tilde{\tau}_i):\tilde{C}
\longrightarrow \pp^1\,\,.$$
We may also write 
\begin{equation}\label{bplusfi}
f=b+\sum_{i=1}^d f_i\,\,,
\end{equation}
where $f_i$ is the morphism $\tilde{C}\to\pp^1$, defined by
$\tilde{h}_i(\tilde{\tau}_i)$, as above.

Note that the base locus of the pencil defined by the sections $1$ and
$b+\sum_i h_i(\tau_i)$ of $\O_C(\sum_iD_i)$ is $Z$ and so
$f:\tilde{C}\to\pp^1$ is the morphism defined by blowing up the locus
of indeterminacy of the rational map~(\ref{ratmap}).

\begin{prop}\label{charf}
Assume that $(C,x)$ is a stable marked curve. Then $(\tilde{C},f)$ is
a stable map of degree $d$. The canonical morphism $p:\tilde{C}\to C$
identifies $(C,x)$ as the stabilization of
$(\tilde{C},\tilde{x})$. The morphism $f:\tilde{C}\to\pp^1$ is
unramified over $\infty\in\pp^1$ and
\begin{equation}\label{prop1}
f^{-1}(\infty)=\sum_{i=1}^d\tilde{D}_i\,\,.
\end{equation}

Assume now that all $\tau_i$ are nowhere vanishing (which implies
that $\tilde{C}=C$) and that all fibers of $C$ are irreducible.
Then we have
\begin{equation}\label{prop2}
\text{$Df({x}_i)({\tau_i})=\frac{\partial}{\partial
z}(\infty)$, for all $i=1,\ldots,d$}\,\,
\end{equation}
and if we let $R$ be the ramification scheme of $f$, then $f\resto
R$ avoids $\infty$, thus is a regular function and we have
\begin{equation}\label{prop3}
b=\frac{1}{2d-2}\tr_{R/T}(f\resto R)\,\,.
\end{equation}
Finally, if $g:C\to\pp^1$ is another family of morphisms of degree $d$
unramified over $\infty$ with Properties~(\ref{prop1}), (\ref{prop2})
and~(\ref{prop3}), then $g=f$.
\end{prop}
\begin{pf}
The stability condition can be checked on fibers of $\tilde{\pi}$. The
only potentially unstable components in such a fiber (over $t$) come
from an exceptional divisor of the blow up.  Such a component
$E_{i,t}$ is isomorphic to $\pp^1$ and intersects a unique
$\tilde{D}_i$ at a unique point which gets mapped to $\infty$. The
component $E_{i,t}$ also intersects $\tilde{V}_{i,C}$, the strict
transform of $V_{i,C}$, in a unique different point, which the map
$f_i$ sends to $0$. Thus $f_i$ is of degree $1$, when restricted to
$E_{i,t}$. Since all the other maps $f_j$, for $j\not=i$, are constant
on $E_{i,t}$, we see that $f\resto E_{i,t}$ is of degree $1$, and
hence that $E_{i,t}$ is a stable component of the map
$(\tilde{C},f)$. 

Since $p:(\tilde{C},\tilde{x})\to(C,x)$ is a morphism of prestable
marked curves, and $(C,x)$ is stable, $p$ has to be the stabilization
morphism.  This follows from the universal property of stabilization
and the fact that every morphism of stable marked curves is an
isomorphism (see \cite{BM}, page~27).

By Proposition~\ref{propdeg1} each of the $f_i$ is unramified over
$\infty$ and maps $\tilde{D}_i$ to $\infty$. Moreover,
$\tilde{\tau}_i$ gets mapped to the canonical tangent vector at
$\infty\in\pp^1$. But $f_j$ for $j\not=i$ is holomorphic (i.e.,
nowhere equal to $\infty$) in a neighborhood of $\tilde{D}_i$, and so
adding it to $f_i$ does not affect these properties of $f_i$ at
$\tilde{D}_i$. Hence $f$ is unramified over $\infty$ and
$f^{-1}(\infty)=\sum_i\tilde{D}_i$. Moreover, the derivative of $f$
has the same behavior at $x_i$ as the derivative of $f_i$, and so
Formula~(\ref{prop2}) follows from Proposition~\ref{propdeg1}.

Now we assume that all $\tau_i$ are nowhere vanishing. This assumption
we make for $\tilde{C}\to C$ to be an isomorphism.  We also assume
that all the fibers of $C$ are irreducible. This has the nice
consequence that, at least locally in $T$, we can find an affine
coordinate $s$ for $C$, such that $s(x_i)$ has no poles, for any $i$ .
More precisely, we can write $C$ as a product $T\times\pp^1$, such
that the sections $x_i$ become functions $x_i:T\to\pp^1$, and we can
arrange things in such a way that $x_i$ avoids $\infty\in\pp^1$, for
all $i$. Then we let $s$ be the affine coordinate for
$\pp^1-\infty$. Now the $s(x_i)$ are regular functions on $T$, which we
abbreviate by 
$$a_i=s(x_i)\,\,.$$
We may now use $s-a_i$ as equation for $D_i$, so that  $h_i$
(see~(\ref{formhi})) becomes
$$h_i=\left(\frac{1}{s-a_i}+ \frac{1}{1-d} \sum_{j\not=i}
\frac{1}{a_j-a_i}\right)ds(a_i)\,\,.$$
We may also use $ds(a_i)$ to trivialize $\omega_i$ and hence $L_i$. We
write
$$\tau_i=q_i\frac{\del}{\del s}(a_i)$$
in this trivialization. Thus 
$$h_i(\tau_i)=q_i\left(\frac{1}{s-a_i}+ \frac{1}{1-d} \sum_{j\not=i}
\frac{1}{a_j-a_i}\right)$$
and hence
\begin{align}
f(s)&=b+\sum_ih_i(\tau_i)\nonumber\\
    &=b+\sum_i\frac{q_i}{s-a_i}+ \frac{1}{1-d} \sum_{i\not=j}
\frac{q_i}{a_j-a_i} \label{ff}
\end{align}
Suppose that $f(s)$ ramifies at
$(r_{\alpha})_{\alpha=1,\ldots,2d-2}$.
Then the $(r_{\alpha})$ are the roots of $f'(s)=0$. A quick
calculation shows that
$$f'(s)=-\sum_i\frac{q_i}{(s-a_i)^2}.$$
and hence $f'(s)=0$ is equivalent to
$$\sum_iq_i\prod_{j\not=i}(s-a_j)^2=0\,\,.$$
Thus we see that
\begin{equation}\label{sr}
\sum_iq_i\prod_{j\not=i}(s-a_j)^2= 
\left(\sum_iq_i\right)\prod_{\alpha=1}^{2d-2}(s-r_{\alpha})\,\,,
\end{equation}
by comparing leading coefficients. By differentiating~(\ref{sr}), we
also get 
\begin{equation}\label{srp}
2\sum_iq_i\sum_{j\not=i}\frac{1}{s-a_j}\prod_{k\not=i}(s-a_k)^2
=\left(\sum_iq_i\right)
\sum_{\alpha}\prod_{\beta\not=\alpha}(s-r_\beta)\,\,.
\end{equation}
Now note that $\sum_{\alpha}\frac{1}{s-r_\alpha}$ is equal to the
right hand side of~(\ref{srp}) divided by the right hand side
of~(\ref{sr}), and so
\begin{equation}\label{sumquot}
\sum_{\alpha}\frac{1}{s-r_\alpha}=
\frac{
2\sum_iq_i\prod_{k\not=i}(s-a_k)^2\sum_{j\not=i}\frac{1}{s-a_j}}{
\sum_iq_i\prod_{k\not=i}(s-a_k)^2}
\end{equation}
We are now ready to compute the trace of $f\resto R$.
\begin{align*}
\sum_{\alpha}f(r_\alpha)&= \sum_\alpha \left(
b+\sum_i\frac{q_i}{r_\alpha-a_i}+ \frac{1}{1-d} \sum_{i\not=j} 
\frac{q_i}{a_j-a_i}\right)\\
&=(2d-2)b-2\sum_{i\not=j} \frac{q_i}{a_j-a_i} - \sum_iq_i\sum_{\alpha}
\frac{1}{a_i-r_\alpha}\\
&=(2d-2)b-2\sum_{i\not=j} \frac{q_i}{a_j-a_i} \\
&\hskip1cm-\sum_iq_i\frac{
2\sum_\ell q_\ell \prod_{k\not=\ell}(a_i-a_k)^2
\sum_{j\not=\ell}\frac{1}{a_i-a_j}}{ 
\sum_\ell q_\ell\prod_{k\not=\ell}(a_i-a_k)^2}\,\,\\
\intertext{(by Equation~(\ref{sumquot}))}
&=(2d-2)b-2\sum_{i\not=j} \frac{q_i}{a_j-a_i} \\
&\hskip1cm-\sum_iq_i\frac{
2 q_i \prod_{k\not=i}(a_i-a_k)^2
\sum_{j\not=i}\frac{1}{a_i-a_j}}{ 
q_i\prod_{k\not=i}(a_i-a_k)^2}\,\,\\
\intertext{(because for $\ell\not=i$ the product vanishes)}
&=(2d-2)b-2\sum_{i\not=j} \frac{q_i}{a_j-a_i} 
-2\sum_iq_i\sum_{j\not=i}\frac{1}{a_i-a_j}\\
&=(2d-2)b\,\,.
\end{align*}
It remains to prove the uniqueness claim.  But it is easy to see that 
$$f(s)=c+\sum_{i=1}^d \frac{q_i}{s-a_i}$$
is the most general degree $d$ morphism which maps $a_i$ to $\infty$,
for all $i$, and satisfies 
$$\left.\frac{d}{ds}\frac{1}{f(s)}\right|_{s=a_i}=\frac{1}{q_i}\,\,,$$
for all $i$.  Thus the  uniqueness follows.
\end{pf}

\subsection{The universal situation}\label{universe}

Consider $\ol{M}_{0,d}$, the scheme of stable curves of genus zero
marked by $\ul{d}=\{1,\ldots,d\}$. Let $\pi:C\to\ol{M}_{0,d}$ be the
universal curve and $x_1,\ldots,x_d$ the universal sections. Let
$\omega_i$ and $L_i$ be line bundles over $\ol{M}_{0,d}$ defined as
above.  Finally, define
$$T=\aaa^1\times\prod_{i=1}^dL_i\,\,,$$
where the product is taken over $\ol{M}_{0,n}$.  Thus $T$ is a
vector bundle of rank $d+1$ over $\ol{M}_{0,d}$.  In particular, $T$
is a smooth scheme of dimension $d-3+d+1=2d-2$. When we pull back any
of the $\ol{M}_{0,d}$-schemes $L_i$ or $C$ to $T$, we endow them with
a subscript $T$ (which we also occasionally omit).

If $S$ is a scheme, then we may think of $S$-valued points of $T$ as
$2d+2$-tuples
$$(C,x,b,\tau)=(C,x_1,\ldots,x_d,b,\tau_1,\ldots,\tau_d)\,\,,$$
where $(C,x)\in\ol{M}_{0,d}(S)$ is a stable marked curve over $S$,
$b\in\aaa^1_S$ is a regular function on $S$ and $\tau_i$ is a tangent
vector of $C$ at $x_i$, or rather a section of $L_i$ over $S$. Similarly,
$S$-valued points of $C_T$ are $2d+3$-tuples
$$(C,x,b,\tau,\Delta),$$
where $(C,x,b,\tau)$ is as above and $\Delta\in C(S)$.

The various projections of $\aaa^1\times\prod_{i=1}^dL_i$ onto its
components define a canonical regular function $b\in\Gamma(T,\O)$ and
canonical sections $\tau_i\in\Gamma(T,L_{i,T})$. These are, of course,
the universal $b$, $\tau$.
As in Section~\ref{degn}, we let $Z=Z_1\cup\ldots\cup Z_d\subset C_T$, where
$Z_i=D_i\cap\{\tau_i=0\}$, and we blow up $C_T$ at $Z$ to obtain
$\tilde{\pi}:\tilde{C}\to T$. 

As in Section~\ref{forevery}, we get for every $i$ a
meromorphic function $f_i=h_i(\tau_i)\in\Gamma\big(C_T,\O(D_i)\big)$, which
we can identify with the meromorphic function
$\tilde{h}_i(\tilde{\tau}_i)\in\Gamma\big(\tilde{C},\O(\tilde{D}_i)\big)$.
Note that $f_i$ is an everywhere non-vanishing section of the line bundle
$\O(\tilde{D}_i)$.
We also get for every $i\not=j$ a regular function
\begin{equation}\label{eta}
\eta_{ij}=f_i(x_j)=h_i(x_j)(\tau_i)=
\tilde{h}_i(\tilde{\tau}_i)(\tilde{x}_j)\,.
\end{equation}
The $\eta_{ij}$ will be very useful in Section~\ref{thevectorfield}.

As in Section~\ref{degn}, we get a meromorphic function
$$b+\sum_{i=1}^d
h_i(\tau_i)\in\Gamma\big(C_T,\O(\textstyle\sum_{i}D_i)\big)\,\,,$$
and an induced morphism
$$f:\tilde{C}\longrightarrow\pp^1\,\,.$$ By Proposition~\ref{charf}
$(\tilde{C},f)$ is a stable map of degree $d$ over $T$, and so we get
an induced morphism
\begin{equation}\label{TtoM}
T\longrightarrow\ol{M}_{0,0}(\pp^1,d)\,\,.
\end{equation}

The symmetric group $S_d$ acts on $\ol{M}_{0,d}$ from the right by 
$$(C,x_1,\ldots,x_d)\cdot\sigma=
(C,x_{\sigma(1)},\ldots,x_{\sigma(d)})\,\,.$$ 
We have compatible actions on $T$ given by 
$$(C,x_1,\ldots,x_d,b,\tau_1,\ldots,\tau_d)\cdot\sigma=
(C,x_{\sigma(1)},\ldots,x_{\sigma(d)},b,
\tau_{\sigma(1)},\ldots,\tau_{\sigma(d)})$$  
and on $C_T$ given by
$$(C,x_1,\ldots,x_d,b,\tau_1,\ldots,\tau_d,\Delta)\cdot\sigma=
(C,x_{\sigma(1)},\ldots,x_{\sigma(d)},b,
\tau_{\sigma(1)},\ldots,\tau_{\sigma(d)},\Delta)\,\,.$$  
The vector bundle $\prod_{i=1}^d L_i$ over $\ol{M}_{0,d}$ descends to
a vector bundle $E=[\prod L_i/S_d]$ of rank $d$ over
$[\ol{M}_{0,d}/S_d]$, which does not split anymore. But still,
$[T/S_d]\to[\ol{M}_{0,d}/S_d]$ is a vector bundle of rank $d+1$.

For all $\sigma\in S_d$, the induced automorphism of $C_T$ identifies
$Z_i$ with $Z_{\sigma(i)}$ and so $Z\subset C_T$ is an invariant
subscheme under the $S_d$-action. Thus we get an induced action of
$S_d$ on the blow-up $\tilde{C}$. This action is compatible with the
projection $\tilde{\pi}:\tilde{C}\to T$ and so we get an induced
prestable curve $[\tilde{C}/S_d]\to[T/S_d]$. Note also that
$f:\tilde{C}\to\pp^1$ is $S_d$-invariant. This follows from
Formula~(\ref{bplusfi}) and the fact that the action of $\sigma$
exchanges $f_i$ and $f_{\sigma(i)}$. Thus we get an induced morphism
$[f/S_d]:[\tilde{C}/S_d]\to\pp^1$, and so we see that (\ref{TtoM})
induces a morphism
\begin{equation}\label{TmodStoM}
[T/S_d]\longrightarrow\ol{M}_{0,0}(\pp^1,d)\,\,.
\end{equation}

\begin{them}\label{th}
The morphism (\ref{TmodStoM}) is an isomorphism onto the open substack
of $\ol{M}_{0,0}(\pp^1,d)$ consisting of stable maps which are
unramified over $\infty\in\pp^1$. 
\end{them}

Before proving the theorem, we prepare a little more. We already
remarked that `unramified over $\infty$' is an open condition on
stable maps to $\pp^1$. So the stable maps unramified over $\infty$
form an open substack $U\subset\ol{M}_{0,0}(\pp^1,d)$. Let $z$ be the
canonical coordinate at $\infty$ on $\pp^1$. If
$(C,f)\in\ol{M}_{0,0}(\pp^1,d)(S)$ is a stable map parametrized by the
scheme $S$, then $f\upst(z)$ defines the closed subscheme
$f^{-1}(\infty)$ of $C$.  

Assume that the stable map $(C,f)$ is unramified over $\infty$. Then
$f^{-1}(\infty)\to S$ is finite \'etale of degree $d$.
We call an isomorphism
$$\phi:\ul{d}\times S\longrightarrow f^{-1}(\infty)$$
an {\em indexing }of $f^{-1}(\infty)$. The local indexings form an
$S$-scheme
$$P=\Isom_S(\ul{d}\times S,f^{-1}(\infty))\,\,,$$ 
which is a principal (right) $S_d$-bundle over $S$. In particular,
$P\to S$ is finite \'etale of degree $d!$.

Let $U'$ be the stack of triples $(C,f,\phi)$, where
$(C,f)\in\ol{M}_{0,0}(\pp^1,d)$ is a stable map unramified over
$\infty$ and $\phi$ is an indexing of $f^{-1}(\infty)$. Thus $U'\to U$
is a principal $S_d$-bundle, in fact, $U'\to U$ is the stack of
indexings of the universal $f^{-1}(\infty)$. It is not difficult to
see that $U'$ is a scheme.

Now consider the stable map $(\tilde{C},f)$ defined over $T$. By
Proposition~\ref{charf} we have that 
$$f^{-1}(\infty)=\sum_{i=1}^{d}\tilde{D}_i\,\,.$$
Thus $(\tilde{C},f)$ comes with a canonical indexing of
$f^{-1}(\infty)$. Therefore, we get a morphism $T\to U'$.  In other
words, the morphism (\ref{TtoM}) lifts in a natural way to
$U'\to\ol{M}_{0,0}(\pp^1,d)$.

Theorem~\ref{th} now follows from the following proposition.

\begin{prop}
The canonical morphism $T\to U'$ is an isomorphism of schemes with
$S_d$-action.
\end{prop}
\begin{pf}
Let us define the inverse of $T\to U'$. Let $(C,f,\phi)$ be an
$S$-valued point of $U'$.  We need to associate to $(C,f,\phi)$ an
$S$-valued point of $T$. Since $U'$ is smooth, we may assume that $S$
is smooth and that the structure morphism $S\to U'$ is \'etale.

Let $S'\subset S$ be the locus over which $f$ does not contract any
components of $C$. Since a contracted component has at least 3 special
points, the complement of $S'$ in $S$ has codimension at least
3. Thus to define a regular function on $S$ is equivalent to defining
a regular function on $S'$ (codimension 2 would suffice for this). Let
$C'\to S'$ be the restriction of $C\to S$ to $S'$.

Let $R$ be the ramification scheme of $f$. Over $S'$ the exact
sequence~(\ref{notses}) is exact on the left also
$$\ses{f\upst\Omega_{\pp^1}}{}{\omega_{C'/S'}}{
}{\omega_{C'/S'}\otimes\O_{R'}}$$ 
and commutes with base change to the fibers of $C'\to S'$. Hence we see
that over $S'$ the pullback $R'=S'\times_SR$ is finite and flat of
degree $2d-2$ over $S'$. 

By assumption, $R'\cap f^{-1}(\infty)=\varnothing$ and so $f\resto R'$
factors through $\aaa^1=\pp^1-\infty$, i.e., $f\resto R'$ is a regular
function on $R'$. Then we can take the trace to get a regular function
$$\tr_{R'/S'}(f\resto R')$$
on $S'$, which extends uniquely to a regular function on $S$, which we
denote by $\tr_{R/S}(f\resto R)$, by abuse of notation. Define
$$b=\frac{1}{2d-2}\tr_{R/S}(f\resto R)\,\,.$$
Thus $b$ is the average of the ramification points of the stable map
$f$.

The isomorphism $\phi:\ul{d}\times S\to f^{-1}(\infty)$ defines $d$
sections $x_1,\ldots,x_d:S\to C$. Consider the derivative
$$Df:\tT_{C/S}\longrightarrow f\upst\tT_{\pp^1}$$
and pull it back via $x_i$ to get an isomorphism
$$Df(x_i):L_i=x_i\upst\tT_{C/X}\longrightarrow x_i\upst
f\upst\tT_{\pp^1}=\tT_{\pp^1}(\infty)\otimes\O_S\,\,,$$
where $\tT_{\pp^1}(\infty)$ is the tangent space of $\pp^1$ at
$\infty$ (see Lemma~\ref{unramlem}).  Taking the preimage of
$\frac{\del}{\del z}(\infty)\tensor 1$ under $Df(x_i)$ yields a
section $\tau_i\in\Gamma(S,L_i)$. 

Stabilizing $(C,x)$ defines a stable marked curve
$$(\ol{C},\ol{x})=(C,x)^\stab\,\,.$$
The stabilization morphism $p:C\to\ol{C}$ induces a homomorphism
$$L_i\to\ol{L}_i=\ol{x}_i\upst\tT_{\ol{C}/S}\,\,,$$
which maps $\tau_i$ to a section $\ol{\tau}_i\in\Gamma(S,\ol{L}_i)$. 

Now we have defined an $S$-valued point
$$(\ol{C},\ol{x}_1,\ldots,\ol{x}_d,b,\ol{\tau}_1,\ldots,\ol{\tau}_d)$$
of $T$, which we declare to by the image of $(C,f,\phi)$, thus
defining a morphism $U'\to T$.

Claim I\@. The composition $T\to U'\to T$ is equal to the identity.

To prove this claim, we start with an $S$-valued point $(C,x,b,\tau)$
of $T$.  We pass to $(\tilde{C},\tilde{x},f)$, which defines a point
of $U'$.  From Proposition~\ref{charf} it follows that $(C,x)$ is the
stabilization of $(\tilde{C},\tilde{x})$, and so the marked curve
associated to $(\tilde{C},\tilde{x},f)$ under $U'\to T$ is $(C,x)$. It
remains to prove that $b=\frac{1}{2d-2}\tr(f\resto R)$ and
$\tau_i=Df(x_i)^{-1}(\frac{\del}{\del z}(\infty))$. But these facts
may be checked on a dense open subset of $S$ (or $T$), and so they
follow from Proposition~\ref{charf}, Equations~(\ref{prop2})
and~(\ref{prop3}).

Claim II\@. The composition $U'\to T\to U'$ is  the identity. 

To prove this claim, we may pass to a dense open substack of $U'$,
because the target $U'$ is separated. 
We start with an $S$-valued point $(C,f,\phi)$ of $U'$,
which defines a prestable marked curve $(C,x)$.  If we assume that $S$
is \'etale over $U'$, then the locus $S'\subset S$ over which $(C,x)$ is
stable is dense in $S$. Over $S'$, neither stabilization nor blowing
up changes $(C,x)$ at all, so that, at least over $S'$, the curve
$(C,x)$ agrees with the one obtained by applying the composition
$U'\to T\to U'$. To check that also the map $f$ does not change after
passing through $U'\to T\to U'$, we may make
$S'$ still smaller, and apply the uniqueness part of
Proposition~\ref{charf}. 
\end{pf}

\begin{rmk}
Denote by $U_a$ the open substack of $\ol{M}_{0,0}(\pp^1,d)$
consisting of maps which are unramified over $a\in\pp^1$.  Obviously,
$U_a$ is isomorphic to $U_\infty$ and hence $U_a\cong T$, for all
$a$. Choosing any $N$ distinct points $a_1,\ldots,a_N\in\pp^1$, where
$N>2d-2$, let $U_i=U_{a_i}$. Then $U_1,\ldots,U_N$ cover
$\ol{M}_{0,0}(\pp^1,d)$. Thus it is possible to obtain
$\ol{M}_{0,0}(\pp^1,d)$ by gluing together $N$ copies of $[T/S_d]$.
Ultimately, this leads to a complete description of the stable map
stack $\ol{M}_{0,0}(\pp^1,d)$ in terms of the stable curve space
$\ol{M}_{0,d}$. 

In this way, many questions about $\ol{M}_{0,0}(\pp^1,d)$ can be
reduced to questions about $\ol{M}_{0,d}$. 
\end{rmk}

\subsubsection{The action of the multiplicative group}

We endow $\pp^1$ with the right action of the group $G$ (see
Remark~\ref{pair}) given in Example~\ref{stupex}.  As usual, we get an
induced (right) action of $B$ on $\ol{M}_{0,0}(\pp^1,d)$ given by
$$(C,f)(a,\lambda)=\big(C,(a,\lambda)\comp f\big)\,.$$
As in Remark~\ref{pair}, this means that we have a $\Gm$-action and an
equivariant vector field on $\ol{M}_{0,0}(\pp^1,d)$.

Since $\infty\in\pp^1$ is a fixed point for the $\Gm$-action on
$\pp^1$, it is obvious that the stack of stable maps unramified over
$\infty$ is invariant under the $\Gm$-action on
$\ol{M}_{0,0}(\pp^1,d)$.  In this section we will determine the
induced $\Gm$-action on $T$.  Since the stack of maps unramified over
$\infty$ is not invariant under all of $B$, we cannot describe an
induced action of $B$ on $T$, but, of course, we can pull back the
equivariant vector field to $T$. Since this is more involved, we shall
postpone it to Section~\ref{thevectorfield}.

\begin{prop}
Let $\Gm$ act on $T$ through scalar multiplication on the vector
bundle $T$ over $\ol{M}_{0,d}$. This $\Gm$-action commutes with the
$S_d$-action, hence induces a $\Gm$-action on $[T/S_d]$. Then the
open immersion~(\ref{TmodStoM}) of Theorem~\ref{th} is
$\Gm$-equivariant.
\end{prop}
\begin{pf}
Let $t\in T$. We have to show that for $\lambda\in\cc^*$ we have
$$(C_t,\lambda\comp f_t)\cong(C_{\lambda t},f_{\lambda t})$$
as stable maps. But to check that two automorphisms of $T$ agree, we
may pass to a dense open subscheme of $T$. So we may assume that
$C_t=C_{\lambda t}=\pp^1$ and that $f_t:\pp^1\to\pp^1$ is given
by~(\ref{ff}). Then the linearity of~(\ref{ff}) in
$(b,q_1,\ldots,q_d)$, which are coordinates on the fibers of
$T\to\ol{M}_{0,d}$, implies that $\lambda\comp f_t=f_{\lambda t}$,
which implies the claim.
\end{pf}

\subsection{Maps to $\pp^n$}\label{mtp}

Let $Y\subset\pp^n$ be the open subvariety defined by
$$Y=\{\langle x_0,\ldots,x_n\rangle\in\pp^n \st\text{$x_0\not=0$ or
$x_1\not=0$}\}\,.$$
The map $\langle x_0,\ldots,x_n\rangle\mapsto\langle x_0,x_1\rangle$
defines a morphism $p:Y\to\pp^1$. 
$$\xymatrix{
Y\dto_p\ar@{^{(}->}[r]& \pp^n\\
\pp^1 & }$$
Let us denote the fiber of $p$ over $\infty$ by $Y_\infty$. Thus
$Y_\infty=\{\langle0,1,y_2,\ldots,y_n\rangle\in\pp^n\}$, and
$Y_\infty$ is canonically identified with $\aaa^{n-1}$.   

We will describe the open substack $U\subset\ol{M}_{0,0}(\pp^n,d)$
consisting of stable maps $(C,f)$ such that $f(C)\subset Y$ and $(C,p\comp
f)$ is unramified over $\infty=\langle 0,1\rangle\in\pp^1$. In fact, we
will show that $U$ is a vector bundle over $\ol{M}_{0,d}$ modulo the action
of $S_d$. 

Let $U_1\subset\ol{M}_{0,0}(\pp^1,d)$ be the open substack of
stable maps unramified over~$\infty$. By definition, we have a cartesian
diagram
$$\xymatrix{
U\ar@{}[dr]|{\square}\dto\ar@{^{(}->}[r] &
\ol{M}_{0,0}(Y,d)\ar@{^{(}->}[r]\dto 
&\ol{M}_{0,0}(\pp^n,d) \\
U_1\ar@{^{(}->}[r] & \ol{M}_{0,0}(\pp^1,d)&}$$

Let 
\begin{equation}\label{strange}
T=\aaa^1\times\left(\prod_{i=1}^d L_i\right)\times
\big(\aaa^1\times\aaa^{d}\big)^{n-1}\,,
\end{equation}
where the product is taken over $\ol{M}_{0,d}$. 
We will write an element of $T(S)$, for a scheme $S$, as
$$(C,x,b,\tau,r)$$
where $(C,x)=(C,x_1,\ldots,x_d)\in\ol{M}_{0,d}(S)$ is a stable marked curve
over $S$, $b=(b_1,\ldots,b_n)$ is an $n$-tuple of regular functions on $S$,
$\tau=(\tau_1,\ldots,\tau_d)$, where $\tau_i$ is a section of $L_i$ over
$S$ and $r=(r_{\nu,i})_{\nu=2,\ldots,n\atop i=1,\ldots,d}$ is a
$d(n-1)$-tuple of regular functions on $S$. The order of the various
coordinates in~(\ref{strange}) is
$$\big(b_1,(\tau_i), b_2, (r_{2,i}), \ldots,
b_n, (r_{n,i})\big)\,.$$
The reason for this order will become clear below.

For notational convenience, we assign the value $1$ to $r_{1,i}$, for
all $i=1,\ldots,d$.

For $\nu=1,\ldots,n$ consider
$$b_\nu+\sum_{i=1}^d
h_i(r_{\nu,i}\tau_i)\in\Gamma\big(C_T,\O(\textstyle\sum_iD_i)\big)\,. $$
Let, as above, $Z=Z_1\cup\ldots\cup Z_d$, where
$Z_i=D_i\cap\{\tau_i=0\}$. After blowing up $Z\subset C_T$ we get a morphism
$$f=\langle1,{\varphi}_1,\ldots,{\varphi}_n\rangle:\tilde{C}\longrightarrow
Y\,,$$ 
where 
\begin{equation}\label{varphi}
\varphi_\nu=b_\nu+\sum_{i=1}^d \tilde{h}_i(r_{\nu,i}\tilde{\tau}_i)
=b_\nu+\sum_{i=1}^dr_{\nu,i}f_i\,.
\end{equation}
This follows easily from the fact that already the sections $1$ and
${\varphi}_1=b_1+\sum_i\tilde{h}_i(\tilde{\tau}_i)=b_1+\sum_if_i$ of
$\O(\sum_i\tilde{D}_i)$ 
generate this invertible 
sheaf, so that any values for ${\varphi}_2,\ldots,{\varphi}_n$ define a
morphism to 
$\pp^n$.  But, for the same reason, this morphism factors through $Y$. 

Thus $(\tilde{C},f)$ is a stable map of degree $d$ parametrized by $T$ and
hence we get a morphism 
\begin{equation}\label{targetn}
T\longrightarrow U\subset\ol{M}_{0,0}(\pp^n,d)\,.
\end{equation}
We will show that (\ref{targetn}) induces an isomorphism $[T/S_d]\cong U$.

Let us write $r_i=\langle 0,1,r_{2i},\ldots,r_{ni}\rangle$. Then we can say
that $r_1,\ldots,r_d\in Y_\infty$ are the points where our stable map
intersects $Y_\infty$ and $\tau_i$ gives the `speed' with which the map
passes through the point $r_i$. The coordinate $b_\nu$ gives the average of
the ramification points of the projection onto the $\nu$-th coordinate
axis.

The next proposition makes this more precise.

\begin{prop}\label{explain}
For any $i=1,\ldots,d$ the composition $f\comp x_i$ defines a morphism 
$$f\comp x_i:T\longrightarrow Y_\infty\,,$$
which is given by $f\comp x_i=r_i$.

Fix a value of $\nu=2,\ldots,n$. Over the locus where
$(r_{\nu,1},\dots,r_{\nu,d})\not=0$ we can compose $f$ with the projection
$p_\nu$ onto the coordinate axis $\{\langle
y_0,0,\ldots,y_\nu,\ldots,0\rangle\}$ to obtain a stable map
$f_\nu=p_\nu\comp f$ unramified over $\infty=\langle
0,\ldots,1,\ldots,0\rangle$. Let $R_\nu\subset\tilde{C}$ be the
ramification scheme of $f_\nu$. Then we have 
$$b_\nu=\frac{1}{2d-2}\tr_{R_\nu/T}(f\resto R_\nu)\,.$$
\end{prop}
\begin{pf}
To check that the two morphism $f\comp x_i$ and $r_i$ agree, we may restrict
to the locus where none of the $\tau_i$ or $r_{\nu,i}$ vanish. Then we need
to prove that 
$$\left.\left(\frac{\varphi_2(s)}{\varphi_1(s)},\ldots,
\frac{\varphi_n(s)}{\varphi_1(s)}\right)\right|_{s=x_i}=
(r_{2,i},\ldots,r_{n,i})\,.$$
This follows from~(\ref{prop2}), using
l'H\^opital's rule. 
The second claim, giving the meaning of the $b_\nu$, follows directly
from~(\ref{prop3}). 
\end{pf}

Let $T_1\to\ol{M}_{0,0}(\pp^1,d)$ be the scheme constructed in
Section~\ref{universe} and called simply $T$ there. By forgetting the new
coordinates we introduced here, we get a morphism $T\to T_1$, which makes
$T$ into a vector bundle over $T_1$. Note that
$$\xymatrix{
T\dto\rto & \ol{M}_{0,0}(Y,d)\dto\\
T_1\rto & \ol{M}_{0,0}(\pp^1,d)}$$
commutes. So letting $\tilde{U}$ be the
fibered product $\tilde{U}=T_1\times_{U_1}U$ we get the diagram
$$\xymatrix{
T\rto\drto&
\tilde{U}\ar@{}[dr]|{\square}\dto\rto & 
U\ar@{}[dr]|{\square}\ar@{^{(}->}[r]\dto&
\ol{M}_{0,0}(Y,d)\dto\ar@{^{(}->}[r]&
\ol{M}_{0,0}(\pp^n,d)\\
&T_1\rto & 
U_1\ar@{^{(}->}[r] & 
\ol{M}_{0,0}(\pp^1,d)&}$$

By Proposition~\ref{VB} the scheme $\tilde{U}$ is also a vector bundle over
$T_1$.

\begin{prop}
The morphism $T\to\tilde{U}$ is an isomorphism of vector bundles over
$T_1$. Hence, $[T/S_d]\cong U$.
\end{prop}
\begin{pf}
To check that $T\to\tilde{U}$ is a morphism of vector bundles is made
easier by the fact that it suffices for this to prove compatibility with
the (linear) $\Gm$-actions. This is proved directly from the definitions of
the two vector bundle structures using Proposition~\ref{explain}. 

Since both vector bundles have the same rank $(d+1)(n-1)$, for the
isomorphism property it suffices to prove strict injectivity, i.e.,
injectivity over every base change to a point of $T_1$. Then the claim also
follows easily from Proposition~\ref{explain}.
\end{pf}

\begin{prop}
Let $\Gm$ act on $T$ by
\begin{multline*}
\big(C,x,b_1,(\tau_i),b_2,(r_{2,i}),\ldots,b_n,(r_{n,i})\big)\lambda\\
=\big(C,x,b_1\lambda,(\tau_i\lambda), b_2\lambda^2,(r_{2,i}\lambda),
\ldots,b_n\lambda^n,(r_{n,i} \lambda^{n-1})\big)\,.
\end{multline*}
In other words, we let $b_\nu$ have weight $\nu$, we let $\tau_i$ have
weight $1$ and we let $r_{\nu,i}$ have weight $\nu-1$. 
Then $T\to\ol{M}_{0,0}(\pp^n,d)$ is $\Gm$-equivariant.
\end{prop}
\begin{pf}
This also follows from Proposition~\ref{explain} using the definition of
the $\Gm$-action on $\pp^n$ from Example~\ref{stupex}.
\end{pf}

\subsection{The degree 2 case}

Since $\ol{M}_{0,2}$ does not exist, the above considerations to not
apply directly to the degree 2 case. We show how to treat this case
here.

Let $U\subset \ol{M}_{0,0}(\pp^n,2)$ be, as above, the open substack
of stable maps $(C,f)$ such that 

(i) $f(C)\subset Y$,

(ii) $p\comp f$ is unramified over $\infty\in\pp^1$ (or $f$ intersects
$Y_\infty$ transversally). 

\noindent We will show how to write $U$ as the stack quotient
$[\aaa^{3n-1}/S_2]$, for a suitable $S_2$-action on $\aaa^{3n-1}$.

Let $T_1=\aaa^2$, with coordinates $b$, $q$. Let $C=T_1\times \pp^1$,
with sections $x_1:T_1\to T_1\times\pp^1; t\mapsto (t,0)$ and
$x_2:T_1\to T_1\times\pp^1;t\mapsto(t,\infty)$. Let $s$ be the
canonical affine coordinate on $\pp^1$, $z=\frac{1}{s}$. Apply the
program of Section~\ref{degn} with $b$, $\tau_1=\frac{\del}{\del
s}(0)$ and $\tau_2=q\frac{\del}{\del z}(\infty)$ to obtain a stable
map of degree 2 from $\tilde{C}$ to $\pp^1$, unramified over $\infty$,
where $\tilde{C}$ is the blow up of $C=T_1\times\pp^1$ at
$Z=Z_2=D_2\cap\{q=0\}$. Note that even though $(C,x_1,x_2)$ is not a
stable curve, $(\tilde{C},x_1,\tilde{x}_2,f)$ is a stable map. Over
$T_1-Z$, the morphism $f$ is given by 
$$f(s)=b+\frac{1}{s}+qs\,.$$

Let $S_2$ act trivially on $T_1$, and on $T_1\times\pp^1$ via the
involution
\begin{align*}
\sigma:T_1\times\pp^1 & \longrightarrow T_1\times\pp^1\\
(b,q,s)&\longmapsto (b,q,\frac{1}{qs})\,.
\end{align*}
Note that $\sigma$ extends uniquely to an automorphism
$\tilde{\sigma}$ of $\tilde{C}$, giving rise to an action of $S_2$ on
$\tilde{C}$. Note also that $f\comp\tilde{\sigma}=f$, and so we get a
stable map
$$\xymatrix{
[\tilde{C}/S_2]\dto\rto^f & \pp^1\\
[T_1/S_2] &}$$
of degree 2. Thus we have a morphism
$[T_1/S_2]\to\ol{M}_{0,0}(\pp^1,2)$. 

\begin{lem}
This morphism $[T_1/S_2]\to\ol{M}_{0,0}(\pp^1,2)$ is an isomorphism
onto the open substack $U$.\qed
\end{lem}

Now let $$T=T_1\times (\aaa^1\times \aaa^2)^{n-1}\,,$$ with additional
coordinates $(b_2,r_{21},r_{22}),\ldots,(b_n,r_{n1},r_{n2})$.
Let $\tilde{C}_T=\tilde{C}\times_{T_1} T$ and define
$\varphi:\tilde{C}_T\to\pp^n$ by $\varphi=\langle
1,f,\varphi_2,\ldots,\varphi_n\rangle$, where 
$$\phi_\nu(s)=b_\nu + r_{\nu 1}\frac{1}{s} + r_{\nu2}\, qs\,,$$
for $\nu=2,\ldots,n$.

Let $S_2$ act on $T$ by fixing $T_1$ and all $b_\nu$, and exchanging
$r_{\nu1}$ with $r_{\nu2}$, for $\nu=2,\ldots,n$. Finally, let $S_2$
act on $\tilde{C}_T$ diagonally, and denote the corresponding
involution of $\tilde{C}_T$ by $\tilde{\sigma}_T$. Then
$\varphi\comp\tilde{\sigma}_T=\varphi$ and so we get an induced stable
map
$\varphi:[\tilde{C}_T/S_2]\to\pp^n$ parametrized by $[T/S_2]$ and
hence a morphism $[T/S_2]\to\ol{M}_{0,0}(\pp^n,2)$.

\begin{them}
The morphism $[T/S_2]\to\ol{M}_{0,0}(\pp^n,2)$ is an isomorphism onto
the open substack $U$. Moreover,
\begin{equation}\label{sq2}
\vcenter{\xymatrix{
[T/S_2]\dto\rto & \ol{M}_{0,0}(Y,2) \dto \\
[T_1/S_2]\rto & \ol{M}_{0,0}(\pp^1,2)}}
\end{equation}
is a pullback diagram of vector bundles.\qed
\end{them}

Moreover, (\ref{sq2}) is $\Gm$-equivariant if we let $\Gm$ act on $T$
by
$$\big(b,q,\ldots,(b_\nu,r_{\nu1},r_{\nu2}),\ldots\big)\cdot\lambda
=\big(b\lambda,q\lambda^2,\ldots,
(b_\nu\lambda^\nu,r_{\nu1}\lambda^{\nu-1},r_{\nu2}\lambda^{\nu-1}),
\ldots\big)\,.$$

\section{The vector field}
\label{thevectorfield}

\subsection{Some deformation theory}

Recall that the group $G=\Ga\rtimes\Gm$ acts on $\pp^n$, as described in
Example~\ref{stupex}.  We get induced actions on the stable map stacks
$\ol{M}_{0,0}(\pp^n,d)$ and $\ol{M}_{0,1}(\pp^n,d)$ by the usual formula
$(C,x,f)\cdot g=(C,x,g\comp f)$. The latter stack we may interpret as the
universal curve $C$ over $\ol{M}_{0,0}(\pp^n,d)$, and so we have a diagram
$$\xymatrix{
C\dto_{\pi}\rto^f & \pp^n\\
{\ol{M}_{0,0}(\pp^n,d)} & }$$
of stacks with $G$-equivariant morphisms.  If we differentiate the
$\Ga$-actions, we obtain $\Gm$-equivariant vector fields $W$ on $\pp^n$,
$V$ on $\ol{M}_{0,0},(\pp^n,d)$ and $U$ on $C$. Of course, $U$ maps to $V$
and $W$ under $D\pi$ and $Df$, respectively. We are concerned with
finding $V$. 

In this section we will show
that $V$ and $U$ are determined uniquely, simply by the fact that
$D\pi(U)=V$ and $Df(U)=W$. The same is then true for the \'etale
$\ol{M}_{0,0}(\pp^n,d)$-scheme $T$: we know that we have found $V\resto T$
if we can find a vector field $U$ on $\tilde{C}$ that simultaneously lifts
$V\resto T $ and $W$. 

\begin{prop}\label{dfiso}
The derivative $Df\colon\tT_C\to f\upst \tT_{\pp^n}$ induces an isomorphism
of vector bundles
$$\pi\lst\tT_C\longrightarrow \pi\lst f\upst\tT_{\pp^n}\,.$$
\end{prop}

We will prove this proposition below.

\begin{cor}\label{herex}
We have a diagram
$$\xymatrix{
&&\Gamma(\pp^n,\tT_{\pp^n})\dto^{f\upst}\\
& \Gamma(\tilde{C},\tT_{\tilde{C}})\rto^{\sim}_{Df}\dto^{D\tilde{\pi}} &
\Gamma(\tilde{C},f\upst\tT_{\pp^n}) \\
\Gamma(T,\tT_{T})\rto^{\sim} &
\Gamma(\tilde{C},\tilde{\pi}\upst\tT_T)\,. &}$$ 
By inverting the isomorphisms we get a homomorphism
\begin{equation*}
\Gamma(\pp^n,\tT_{\pp^n})\longrightarrow\Gamma(T,\tT_T)\,,
\end{equation*}
which maps $W$ to $V$.\qed
\end{cor}

To prove Proposition~\ref{dfiso} we start by recalling a few general facts
about tangent complexes.

\begin{lem}\label{tcc}
Let 
$$\xymatrix{
X\dto_p\rto^f & Y\dto^g\\
W\rto^q & Z}$$
be a cartesian diagram of smooth stacks. Let $r\cong q\comp p\cong
g\comp f$. Then the diagram 
$$\xymatrix@=.5cm{
r\upst\tT_Z \rrto^{K.S.(g)}\ddto_{K.S.(q)} &&
{f\upst\tT\com_{Y/Z}[1]}\ar@{=}[r]&{\tT\com_{X/W}[1]}\ar[dd]\\ \\
p\upst\tT\com_{W/Z}[1]\ar@{=}[r]&{\tT\com_{X/Y}[1]}\ar[rr]
&&\tT\com_X[1]&}$$
anti-commutes in the derived category of $\O_X$-modules. Here `$K.S.$'
stands for `Kodaira-Spencer' map.
\end{lem}
\begin{pf}
Note that the composition
$$\xymatrix{r\upst\tT_Z\rrto^{K.S.(r)}& & \tT_{X/Z}[1]\rto^\alpha &
\tT_X[1]}$$ is zero. On the other hand, we have
$$\tT_{X/Z}[1]=\tT_{X/Y}[1]\oplus\tT_{X/W}[1]$$
and under this decomposition we have $K.S.(r)=K.S.(q)\oplus
K.S.(g)$. Applying $\alpha$, we get
$$0=\alpha\comp K.S.(r)=\alpha\comp K.S.(q)+\alpha\comp K.S.(g)\,.$$
This is what we wanted to prove.
\end{pf}

\begin{lem}\label{thislem}
Again, considering a cartesian diagram as in the previous lemma, there
is a canonical homomorphism of distinguished triangles
$$\xymatrix{
{f\upst\tT\com_{Y/Z}}\rto\dto & {f\upst \tT\com_Y}\rto\ar@{=}[d] & {f\upst
g\upst\tT\com_Z}\rto^{K.S.(g)}\dto_{K.S.(q)} &
{f\upst\tT\com_{Y/Z}[1]}\dto \\ 
{\tT\com_X}\rto^\beta & {f\upst\tT\com_Y} \rto &
{\tT\com_{X/Y}[1]}\rto^{-\beta[1]}  &
{\tT\com_X[1]}\,.}$$
The upper triangle is the pullback under $f$ of the distinguished
triangle associated with $g$, the lower triangle is a shift of the
distinguished triangle associated with $f$.
\end{lem}
\begin{pf}
We have to show that the three squares commute. The last one is
Lemma~\ref{tcc}. 

In any cartesian diagram as the one under consideration, there is a
homomorphism from the distinguished triangle for $p$ to the pullback
under $f$ of the distinguished triangle for $g$. Similarly, for $f$
and $q$.  Picking out appropriate commutative squares from these
homomorphisms of distinguished triangles proves the commutativity of
the other two squares.
\end{pf}

Consider the morphism of stacks 
$$h:\ol{M}_{0,0}(\pp^n,d)\longrightarrow\MM_{0,0}\,.$$
Here $\MM_{0,0}$ is the Artin stack of prestable curves of genus
zero. The morphism $h$ is given by forgetting the map, retaining only
the prestable curve (and not stabilizing). The fiber of $h$ over a
prestable curve $C$ is equal to an open subscheme of $\Mor(C,\pp^n)$,
the scheme of  morphisms from $C$ to $\pp^n$.

Deformation theory for morphisms shows that $\Mor(C,\pp^n)$ is smooth
with tangent space $H^0(C,f\upst \tT_{\pp^n})$ at the point
$f:C\to\pp^n$ of $\Mor(C,\pp^n)$. It follows that $h$ is smooth with
relative tangent bundle
$$\tT_{\ol{M}/\MM}=\pi\lst f\upst\tT_{\pp^n}\,.$$

We shall now apply Lemma~\ref{thislem} to the cartesian square
$$\xymatrix{
C\dto\rrto^\pi && \ol{M}_{0,0}(\pp^n,d)\dto^h\\
\CC\rrto^{\tilde{\pi}} && \MM_{0,0}}$$
where $\CC\to\MM_{0,0}$ is the universal curve over $\MM_{0,0}$. 
We obtain the homomorphism of distinguished triangles
$$\xymatrix{
{\pi\upst\tT_{\ol{M}/\MM}}\rto\dto & {\pi\upst
\tT_{\ol{M}}}\rto\dto & {\pi\upst 
h\upst\tT\com_{\MM}}\rto\dto & {\pi\upst\tT_{\ol{M}/\MM}[1]}\dto \\
{\tT_C}\rto & {\pi\upst\tT_{\ol{M}}} \rto & {\tT\com_{C/\ol{M}}[1]}\rto &
{\tT_C[1]}}$$
and by adjointness, the homomorphism
$$\xymatrix{
{\tT_{\ol{M}/\MM}}\rto\dto^{\phi} & {
\tT_{\ol{M}}}\rto\dto^{\phi'} & { 
h\upst\tT\com_{\MM}}\rto\dto^{\phi''} & {\tT_{\ol{M}/\MM}[1]}\dto \\
{R\pi\lst\tT_C}\rto & {R\pi\lst\pi\upst\tT_{\ol{M}}} \rto &
{R\pi\lst\tT\com_{C/\ol{M}}[1]}\rto & 
{R\pi\lst\tT_C[1]}\,.}$$
By the projection formula
we have $R\pi\lst\pi\upst\tT_{\ol{M}}=\tT_{\ol{M}}\otimes
R\pi\lst\O_{C}=\tT_{\ol{M}}$. Therefore $\phi'$ is an
isomorphism. We also have that $R\pi\lst\tT\com_{C/\ol{M}}[1]=h\upst
R\tilde{\pi}\lst\tT\com_{\CC/\MM}[1]$. Thus $\phi''$ may be viewed as the
pullback under $h$ of the Kodaira-Spencer homomorphism $\tT\com_{\MM}\to
R\tilde{\pi}\lst\tT\com_{\CC/\MM}[1]$, which is an isomorphism by the
deformation theory of prestable curves.

Since both $\phi'$ and $\phi''$ are isomorphisms, so is
$\phi:\tT_{\ol{M}/\MM}\to R\pi\lst\tT_C$.

Finally, let us consider the composition of homomorphisms of vector
bundles on $C$
$$\pi\upst\tT_{\ol{M}/\MM}=
\tT_{C/\CC}\longrightarrow\tT_C\stackrel{Df}{\longrightarrow} 
f\upst\tT_{\pp^n}\,.$$
Again, by adjointness, we get an induced homomorphism in the  derived
category of $\ol{M}_{0,0}(\pp^n,d)$
$$\tT_{\ol{M}/\MM}\stackrel{\phi}{\longrightarrow}
R\pi\lst\tT_C\stackrel{\psi}{\longrightarrow} 
R\pi\lst f\upst\tT_{\pp^n}\,.$$
We just saw that $\phi$ is an isomorphism, and the composition
$\psi\comp\phi$ is an isomorphism by the deformation theory of
morphisms. Therefore, $\psi$ is an isomorphism. Since $R\pi\lst
f\upst\tT_{\pp^n} =\pi\lst f\upst\tT_{\pp^n}$, we deduce that
$R\pi\lst\tT_C=\pi\lst\tT_C$, and that we have an isomorphism of
vector bundles 
$$\pi\lst\tT_C\longrightarrow\pi\lst f\upst\tT_{\pp^n}\,$$
which is induced by the derivative $Df:\tT_C\to f\upst\tT_{\pp^n}$. 
This finishes the proof of Proposition~\ref{dfiso}.

\subsection{The degree 3 case}

We will now use Corollary~\ref{herex} to determine the vector field
$V$ on $T$ in the case of $d=3$.

So let $d=3$. Then $\ol{M}_{0,3}=\spec\cc$. The universal curve $C$ over
$\ol{M}_{0,3}$ is isomorphic to $\pp^1$, with 3 marked points
$x_1,x_2,x_3\in\pp^1$. The vector bundle
$T$ over $\ol{M}_{0,3}=\spec\cc$ is just a vector space. We have
canonically
$$T=\big(\cc\oplus L_1 \oplus L_2 \oplus
L_3\big)\oplus\big(\cc\oplus\cc^3\big)^{(n-1)},$$ 
where $L_i=\tT_{\pp^1}(x_i)$ is the tangent space of $\pp^1$ at a marked
point.

Every vector space $E$ has a canonical vector field on it, namely the
vector field that takes the value $e\in E$ at the point $e\in E$. This
makes sense, because every tangent space of $E$ is canonically identified
with $E$. If we choose a basis $(e_i)$ for $E$, then this vector field is
given by $\sum_i x_i\frac{\del}{\del x_i}=\sum_i x_i e_i$, where $(x_i)$
are the coordinates on $E$ induced by the basis $(e_i)$. We denote the
canonical vector field on $L_i$ by $\tau_i$, for $i=1,2,3$. Because of the
product decomposition of $T$, we can also think of the $\tau_i$ as
canonical vector fields on $T$. 

Recall that the first coordinate of every quadruple of coordinates on $T$
is called
$b_\nu$, where ${\nu=1,\ldots,n}$.  The latter three coordinates in the
latter $n-1$ quadruples of coordinates are called
$(r_{\nu,i})_{\nu=2,\ldots,n,\, i=1,2,3}$. Recall also that we write
$r_{1,i}=1$, $i=1,2,3$, for notational convenience.  For further notational
convenience we write $r_{\nu,0}=b_{\nu}$, for $\nu=1,\ldots,n$.  This then
defines $r_{\nu,i}$ for all $\nu=1,\ldots,n$ and all $i=0,1,2,3$.  For every
$\nu=1,\ldots,n$ we combine $r_{\nu,0}\ldots,r_{\nu,3}$ into the column
vector $r_\nu$. Thus
$$r_1=\begin{pmatrix}b_1\\1\\1\\1\end{pmatrix}\quad\text{and}\quad
r_\nu=\begin{pmatrix}b_{\nu}\\r_{\nu,1}\\r_{\nu,2}\\r_{\nu,3}\end{pmatrix}
\quad\text{for $\nu=2,\ldots,n$}$$

Because of the product decomposition of $T$, every vector field on $T$ is
in a canonical way a sum of $4n$ components. All but three are canonically
identified with regular 
functions on $T$. We combine groups of four components into one as follows.
$$V_1=V_{1,0}\frac{\del}{\del b_1}+\sum_{i=1}^3 V_{1,i}\,,$$
and for $\nu=2,\ldots,n$
$$V_\nu=V_{\nu,0}\frac{\del}{\del b_\nu}+\sum_{i=1}^3
V_{\nu,i}\frac{\del}{\del r_{\nu,i}}\,.$$ 

It will be convenient to introduce the following abbreviations:
$$\eta_i^{(2)}=-\prod_{j\not=i}\eta_{ij}$$
and for all $\mu\geq1$
$$\theta_i^{(\mu)}=\sum_{j\not=i}\eta_{ij}^{\mu-1}\eta_{ji}\,,$$ 
(so that the upper indices denote the degree). Here the $\eta_{ij}$ are the
canonical regular functions on $T$ introduced in~(\ref{eta}). 
Finally, let 
$$\ell_i=r_{2,i}-2b_1-4\theta_i^{(1)}\,,$$
for $i=1,2,3$.

\begin{prop}\label{vfe}
With this notation, the vector field $V$ on $T$ is
given by
\begin{align*}
V_1&=r_2-Er_1+
\begin{pmatrix}0\\ \ell_1\tau_1\\ \ell_2\tau_2\\
\ell_3\tau_3\end{pmatrix}=
\begin{pmatrix}b_2-b_1^2+\sum_i\eta_i^{(2)}+5\sum_i\theta_i^{(2)}\\
\ell_1\tau_1\\ \ell_2\tau_2\\ \ell_3\tau_3\end{pmatrix}\,,\\
V_\nu&=r_{\nu+1}-Er_\nu\,,\quad\text{for $\nu=2,\ldots,n-1$}\,,\\
V_n&=-Er_n\,.
\end{align*}
Here $E$ is the $4\times4$ matrix 
$$E=\left(\begin{array}{cccc}
b_1 & \eta_1^{(2)}+5\theta_1^{(2)} &
\eta_2^{(2)}+5\theta_2^{(2)} &
\eta_3^{(2)}+5\theta_3^{(2)} \\
1 & \ell_i + b_1 +3\theta_1^{(1)} & \eta_{21} & \eta_{31} \\
1 & \eta_{12}& \ell_2+b_1+3\theta_2^{(1)}  & \eta_{32} \\
1 & \eta_{13} & \eta_{23} & \ell_3+b_1+3\theta_3^{(1)}
\end{array}\right)\,.$$
In the case $n=1$, the formula for $V_1$ is used, not the formula for
$V_n$. Moreover, $b_2$ and $r_{2,i}$ are set equal to $0$ in this case.
\end{prop}
\begin{pf}
It suffices to check these formulas on a dense open subset of $T$. So we
shall restrict to the subscheme of $T$ where none of the $\tau_i$ vanish,
so that we may assume that $\tilde{C}=C_T=T\times\pp^1$. Then all we have to
do is exhibit a vertical vector field $\ol{U}$ for the projection
$T\times\pp^1\to T$ such that $Df(V+\ol{U})=W$. Then the vector field $U$
on $\tilde{C}$ is given by $U=V+\ol{U}$, over this open subscheme of
$\tilde{C}$. 

The relative vector field $\ol{U}$ is very easy to describe. For every
point $t$ of $T$, we have given 3 tangent vectors on $\pp^1$, namely the 
$\tau_i(t)\in\tT_{\pp^1}(x_i)$. On $\pp^1$ there is a unique vector field
$\ol{U}(t)$ taking this prescribed value $\tau_i(t)$ at the point $x_i$,
for $i=1,2,3$.  This defines $\ol{U}$. 

The details follow below.
\end{pf}

\begin{lem}\label{olu}
Let $T$ be a scheme and $\tau_i\in\Gamma(T,L_i)$, for $i=1,2,3$ a section
of $L_i$ over $T$. Then the formula
$$\ol{U}=\sum_{i=1}^3\left(\prod_{j\not=i}
\left(1-\frac{s_i}{s_i(x_j)}\right)\right)\tau_ids_i(x_i)\frac{\del}{\del
s_i}\,,$$ 
defines a relative vector field $\ol{U}\in\Gamma(T\times
\pp^1,\tT_{T\times\pp^1/T})=\Gamma(T\times\pp^1,\tT_{\pp^1})$, which has
the property 
that $\ol{U}(x_i)=\tau_i$, for $i=1,2,3$. Here $s_i$ is a degree one
meromorphic function on $\pp^1$, which vanishes at $x_i$ (i.e., a parameter
at $x_i$).

Let $f_i=h_i(\tau_i)$, for $i=1,2,3$ be the meromorphic function on
$T\times\pp^1$ considered above (see Corollary~\ref{hi}).
Then we have
$$\ol{U}(f_i)=2\theta_i^{(2)}+\eta_i^{(2)}
+2\theta_i^{(1)}f_i-f_i^2\,,$$
for every $i=1,2,3$.
\end{lem}
\begin{pf}
It is not difficult to see that the formula for $\ol{U}$ defines a vector
field on $\pp^1$, no matter the choice of the $s_i$, even if
$s_i(x_j)=\infty$, for some $j\not=i$. The verification that
$\ol{U}(x_i)=\tau_i$ is then easy.

The calculation of $\ol{U}(f_i)$ is somewhat tedious. It can be done by
choosing a coordinate $s$ for $\pp^1$, which does not take the value
$\infty$ at any of the $x_i$. Then write $a_i=s(x_i)$ and
$s_i=s-a_i$. Finally, write $q_i=\tau_ids_i(x_i)$. With these choices we
have 
$$\ol{U}=\sum_{i=1}^3
q_i\left(\prod_{j\not=i}\frac{s-a_j}{a_i-a_j}\right)\frac{\del}{\del
s}$$
and
$$f_i=q_i\left(\frac{1}{s-a_i}-
\frac{1}{2}\sum_{j\not=i}\frac{1}{a_j-a_i}\right)\,.$$
Thus
$$\ol{U}(f_i)=-\sum_{k=1}^3
q_k\left(\prod_{j\not=k}\frac{s-a_j}{a_k-a_j}\right)
\frac{q_i}{(s-a_i)^2}$$
and this can be compared to the given formula for $\ol{U}(f_i)$.
\end{pf}

\begin{lem}\label{fifj}
With the same notation, we have for $i\not=j$
$$f_if_j=\eta_{ji}f_i+\eta_{ij}f_j+3\eta_{ij}\eta_{ji}\,.$$
In particular, summing over all $i$ and $j$ such that $i\not=j$,
$$\sum_{i\not=j}f_i f_j=
2\sum_{i\not=j}\eta_{ji}f_i+3\sum_{i\not=j}\eta_{ij}\eta_{ji}\,.$$
\end{lem}
\begin{pf}
This can be proved similarly.
\end{pf}

In view of the formula for $W$ given in Example~\ref{stupex}, the following
corollary now finishes the proof of Proposition~\ref{vfe}.  

\begin{cor}\label{silcor}
If we were to define a vector field $V$ on $T$ by the formulas of
Proposition~\ref{vfe}, then we would have
\begin{align*}
V(\varphi_\nu)+\ol{U}(\varphi_\nu)&=\varphi_{\nu+1}-\varphi_1\varphi_\nu\,,\\
\intertext{for $\nu<n$ and}
V(\varphi_n)+\ol{U}(\varphi_n)&=-\varphi_1\varphi_n\,.
\end{align*}
Here the $\varphi_\nu$ are the components of the morphism
$f:\tilde{C}\to\pp^n$ given as in (\ref{varphi}) by
$\varphi_\nu=b_\nu+\sum_{i=1}^3r_{\nu,i}f_i$.
\end{cor}
\begin{pf}
We note
that $\tau_j(f_i)=\delta_{ij}f_i$ (Kronecker delta). Moreover,
$V(f_i)=\ell_if_i$. Therefore, we need to
prove the 
following formulas. First
$$
V_{1,0}+\sum_i\ell_if_i+\sum_i\ol{U}(f_i)=
b_2+\sum_ir_{2,i}f_i-\big(b_1+\sum_if_i\big)^2 \,,$$
then for every $1<\nu<n$
\begin{multline*}
V_{\nu,0}+\sum_ir_{\nu,i}\ell_if_i+\sum_iV_{\nu,i}f_i+\sum_ir_{\nu,i}
\ol{U}(f_i)\\= 
b_{\nu+1}+\sum_ir_{\nu+1,i}f_i-
\big(b_1+\sum_if_i\big)\big(b_\nu+\sum_ir_{\nu,i}f_i\big) \,,
\end{multline*}
and finally,
\begin{multline*}
V_{n,0}+\sum_ir_{n,i}\ell_if_i+\sum_iV_{n,i}f_i+\sum_ir_{n,i}
\ol{U}(f_i)\\= 
-\big(b_1+\sum_if_i\big)\big(b_n+\sum_ir_{n,i}f_i\big) \,.
\end{multline*}
All of these formulas follow easily from Lemmas~\ref{olu}
and~\ref{fifj}. 
\end{pf}

To make the relationship between vector fields and regular functions on the
$L_i$ explicit, we identify the universal curve $C$ with $\pp^1$, with affine
coordinate $s$. Then we choose for $x_1$, $x_2$ and $x_3$, the points $0$,
$1$ and $\infty$.  We choose local parameters $\frac{1}{2}s$ at $0$,
$\frac{s-1}{2s}$ 
at $1$ and $\frac{1}{2-2s}$ at $\infty$. This trivializes $L_1$, $L_2$ and
$L_3$, and we write, as usual, the induced  coordinates on $L_i$ as $q_i$. 
With this notation the universal map is given by
$$\varphi_\nu(s)=
b_\nu+r_{\nu,1}\,q_1\frac{2-s}{s}+ r_{\nu,2}\,q_2\frac{s+1}{s-1} +
r_{\nu,3}\,q_3(1-2s)\,,$$
for all $\nu=1,\ldots,n$. Moreover, we have $\tau_i=q_i$ and 
\begin{align*}
\eta_{12}&=q_1 & 
\eta_{13}&=-q_1\\
\eta_{21}&=-q_2 & 
\eta_{23}&=q_2\\
\eta_{31}&=q_3 &
\eta_{32}&=-q_3\,.
\end{align*}
The matrix $E$ then takes the following shape:
$$E=\left(\begin{array}{cccc}
b & q_1\big(q_1-{5}(q_2+q_3)\big) &
q_2\big(q_2-{5}(q_1+q_3)\big) &
q_3\big(q_3-{5}(q_1+q_2)\big) \\
1 & \ell_1+b+3(q_3-q_2) & -q_2 & q_3 \\
1 & q_1& \ell_2+b+3(q_1-q_3)  & -q_3 \\
1 & -q_1 & q_2 & \ell_3+b+3(q_2-q_1) 
\end{array}\right)\,,$$
where we have dropped the index on $b_1$.

\begin{cor}\label{tring}
We have
$$\hh^0(T,K_V\com)=\cc[b,q_1,q_2,q_3,\ell_1,\ell_2,\ell_3]/
(q_1\ell_1,q_2\ell_2,q_3\ell_3,R_{n+1})\,$$
where $R_{n+1}$ denotes the quadruple of relations given by the matrix
equation 
$$R_{n+1}=E^n r_1\,.$$
The degree of all generators $b,q_i,\ell_i$ is one. The degrees of the
components of $R_{n+1}$ are $n+1$, $n$, $n$ and $n$, respectively. The
group $S_3$ acts 
as follows: $\sigma b= b$, $\sigma q_i =\sign(\sigma)q_{\sigma(i)}$,
$\sigma \ell_i=\ell_{\sigma(i)}$, for $\sigma\in S_3$. Note that $S_3$ acts
by a similar pattern on the seven relations.
\end{cor}
\begin{pf}
We have 
$$\hh^0(T,K_V\com)= \cc[b,(q_i)_{i=1,2,3},(r_{\nu,i})_{\nu\geq2,i\geq0}]/
(V_{\nu,i})\,.$$ 
The first four relations, $V_1=0$, give the relations
$(q_i\ell_i)_{i=1,2,3}$ and $r_2=E r_1$. The latter is only one relation
and can be used to eliminate $b_2$, but not the $r_{2,i}$, for
$i=1,2,3$. Then the other relations $V_\nu=0$, for $2\leq\nu<n$, recursively
eliminate all $r_{\nu,i}$, for $\nu\geq3$. The last quadruple of relations
can then be expressed as $R_{n+1}=0$.  Finally, we have replaced the
generators $r_{2,i}$ by $\ell_i$, for $i=1,2,3$.
\end{pf}

\subsection{Chern classes}\label{chern}

Recall the diagram
$$\xymatrix{
C\dto_\pi\rto^f & \pp^n\\
\ol{M}_{0,0}(\pp^n,d) &}$$
involving the universal curve and the universal map. For every $m>0$ we
consider the vector bundle
$$E_m=\pi\lst f\upst\O(m)\,$$
on the stable map stack $\ol{M}_{0,0}(\pp^n,d)$. Our goal in this section
is to show how to determine the characteristic classes
$c_{\tilde{V}}^Q\in\hh^0(T,K_V\com)$ of $E_m$ (see Section~\ref{Chern}).

We restrict ourself to the case that $d=3$. First, let us show
that $E_m$ is trivial over $T$. 

\begin{lem}\label{basis}
Let $f_i=\tilde{h}_i(\tilde{\tau}_i)$, $i=1,2,3$,  be the canonical sections
of 
$\O_{\tilde{C}}(\sum_i\tilde{D}_i)$ constructed in Section~\ref{SecPar}. Via
the canonical identification
$f\upst\O_{\pp^n}(1)=\O_{\tilde{C}}(\sum_i\tilde{D}_i)$ (see the definition
of $f$ in Section~\ref{mtp}) we consider the $f_i$ as global sections of
$\pi\lst f\upst\O(1)\resto T$. Then the family $(1,f_1,f_2,f_3)$
trivializes this bundle.  Moreover, the family
$(1,f_1,f_2,f_3,\ldots,f_1^m,f_2^m,f_3^m)$ of global sections of 
$\O_{\tilde{C}}(\sum_im\tilde{D}_i)$ trivializes the bundle $E_n\resto
T$. Finally, $(1,f_1,f_2,f_3,\ldots,f_1^m,f_2^m,f_3^m)$ is a
$\Gm$-eigenbasis and we have ${^\lambda}f_i^\nu=\lambda^\nu f_i$, for
all $i=0,\ldots,m$.

\end{lem}
\begin{pf}
We have to show that for every geometric point $t$ of $T$ the family
$(1,f_1,f_2,f_3)$  induces a
basis of the fiber $\big(\pi\lst
f\upst\O(1)\big)(t)=H^0(\tilde{C}_t,\O(\sum_i\tilde{x}_i)$. But this is
true, because $f_i(t)$ has a pole of order one at $\tilde{x}_i$, by
construction. The `moreover' follows similarly.
\end{pf}

Recall that the vector field $W$ on $\pp^n$ acts on $\O(m)$ (see
Examples~\ref{nablaex} and~\ref{tildeW}). Via Lemmas~\ref{tildepull}
and~\ref{tildepush} we get an induced action $\tilde{V}$ of $V$ on $E_m$,
for all $m$. 
Of course, all these actions $\tilde{V}$ are $\Gm$-equivariant. To find
the characteristic class $c_{\tilde{V}}^Q(E_m)\resto T$ (where $Q$ is
an invariant 
polynomial function on $(3m+1)\times(3m+1)$-matrices), we only need to
compute the matrix of $\tilde{V}$ with respect to the basis of $E_m\resto
T$ given by Lemma~\ref{basis}, and then evaluate $Q$ on this matrix. This
matrix can be computed from the following lemma, using Lemmas~\ref{olu}
and~\ref{fifj}. 

\begin{lem}
Let $\tilde{V}$ be the action of $V$ on $E_m$. Then we have, for all
$\mu=0,\ldots,m$ 
$$\tilde{V}(f_i^\mu)=(mb+\mu\ell_i)
f_i^{\mu}+ \mu f_i^{\mu-1}\ol{U}(f_i)+
mf_i^{\mu+1}+m\sum_{m\not=i}f_i^\mu f_j\,.$$
\end{lem}
\begin{pf}
It suffices to check this formula on a dense open subset of $T$. So we may
assume that $\tilde{C}=C=T\times\pp^1$, and that the vector field $U$ on
$\tilde{C}$ splits as $U=V+\ol{U}$, where $\ol{U}$ is the relative vector
field defined in Lemma~\ref{olu}. Denote by $\tilde{U}$ the action of $U$
on $f\upst\O(m)$. Then we have, by Lemmas~\ref{tildepull}
and~\ref{tildepush} 
\begin{align*}
\tilde{V}(f_i^\mu) &=
\tilde{U}(f_i^\mu) \\
&= \tilde{U}(f_i^\mu 1)\\
&= U(f_i^\mu)1+f_i^\mu\tilde{U}(1)\\
&= \mu f_i^{\mu-1}U(f_i)+f_i^\mu\tilde{U}(f\upst x_0^m)\\
&= \mu f_i^{\mu-1}V(f_i)+\mu f_i^{\mu-1}\ol{U}(f_i)+f_i^\mu
f\upst\big(\tilde{W}(x_0^m)\big)\\ 
&= \mu f_i^{\mu}\ell_i+\mu f_i^{\mu-1}\ol{U}(f_i)+f_i^\mu
f\upst\big(mx_0^{m-1}\tilde{W}(x_0)\big)
\intertext{because, as we noted in the proof of Corollary~\ref{silcor}, we
have $V(f_i)=\ell_if_i$,}
&= \mu \ell_i f_i^{\mu}+\mu f_i^{\mu-1}\ol{U}(f_i)+f_i^\mu
f\upst(mx_0^{m-1}x_1)\\ 
&= \mu \ell_i f_i^{\mu}+\mu f_i^{\mu-1}\ol{U}(f_i)+mf_i^\mu
f\upst(x_1)\\ 
&= \mu\ell_i f_i^{\mu}+\mu f_i^{\mu-1}\ol{U}(f_i)+mf_i^\mu
\big(b+\sum_j f_j\big)\\
&= \mu\ell_i f_i^{\mu}+\mu f_i^{\mu-1}\ol{U}(f_i)+mbf_i^\mu
+mf_i^{\mu+1}+m\sum_{j\not=i}f_i^\mu f_j\,,
\end{align*}
which is what we wanted to prove.
\end{pf}

\begin{cor}\label{matrix}
The matrix of $\tilde{V}$ acting on $E_m$, with respect to the basis
$(1,f_1,f_2,f_3,\ldots,f_1^m,f_2^m,f_3^m)$ is determined by
\begin{align*}
\tilde{V}(1)&=mb+m(f_1+f_2+f_3)\,,\\
\tilde{V}(f_i)&=(3m+2)\theta_i^{(2)}+\eta_i^{(2)}\\
&+m\sum_{j\not=i}\eta_{ij}f_j
+\big(\ell_i+mb+(m+2)\theta_i^{(1)}\big)f_i + (m-1) f_i^2\,,
\end{align*}
and for all $\mu\geq2$,
\begin{multline*}
\tilde{V}(f_i^\mu)=3m\theta_i^{(\mu+1)}+m\sum_{j\not=i}\eta_{ij}^\mu f_j 
+4m\sum_{\alpha=1}^{\mu-2}\theta_i^{(\mu-\alpha+1)}f_i^\alpha\\
+\big((2\mu +4m)\theta_i^{(2)}+\mu\eta_i^{(2)}\big)f_i^{\mu-1}
+\big(\mu \ell_i+mb+(m+2\mu)\theta_i^{(1)}\big)f_i^\mu + (m-\mu)f_i^{\mu+1}\,,
\end{multline*}
\end{cor}
\begin{pf}
We just need to plug in the formula for $\ol{U}$ given in Lemma~\ref{olu}
and the formula
$$f_i^\mu f_j=\eta_{ji}
f_i^\mu+4\sum_{\alpha=1}^{\mu-1}\eta_{ij}^{\mu-\alpha} f_i^\alpha
+\eta_{ij}^\mu f_j+3\eta_{ij}^\mu\eta_{ji}\,,$$
which holds for $i\not=j$ and $\mu\geq1$, and is proved by
induction from Lemma~\ref{fifj}.
\end{pf}

\begin{rmk}
Note that the matrix for $\tilde{V}$ on $E_1$ is the same as the matrix $E$,
above, determining $V$ itself. This is due to the following phenomenon.  We
may interpret the vector bundle $E_1$ as an open substack of
$\ol{M}_{0,0}(\pp^{n+1},d)$. That gives us two vector fields on $E_1$. One
is the reinterpretation of the action $\tilde{V}$ as a vector field on
$E_1$ (see Remark~\ref{rein}), the other is the vector field $V$ on
$\ol{M}_{0,0}(\pp^{n+1},d)$, restricted to $E_1$.  These two vector fields
are almost, but not quite, equal.
\end{rmk}

\begin{cor}\label{first}
We have
$$c_1(E_m)=m(3m+1)b+\frac{1}{2}m(m+1){\textstyle\sum}\ell_i$$
\end{cor}
\begin{pf}
The first Chern class of $E_m$ is the trace of the matrix
of $\tilde{V}$ on $E_m$. Also note that $\sum_i\theta_i^{(1)}=0$.
\end{pf}

\begin{cor}\label{second}
Let $M_m$ be the matrix of $\tilde{V}$ described in Corollary~\ref{matrix}.
Then we have
$$c_2(E_m)=\frac{1}{2}\big(c_1^2(E_m)-\tr(M_m^2)\big)$$
and
\begin{multline*}
\tr(M_n^2)=m^2(3m+1)b^2 + m^2(m+1)b{\textstyle\sum}\ell_i + 
\frac{1}{6}m(m+1)(2m+1){\textstyle\sum}\ell_i^2 \\
+m(3m+1)(3m+2)\big({\textstyle\sum}\eta_i^{(2)}+
{\textstyle\sum}\theta_i^{(2)}\big) 
+\frac{1}{3}m(m+1)(7m+2){\textstyle\sum}\ell_i\theta_i^{(1)}\,.
\end{multline*}
\end{cor}
\begin{pf}
We use the formula
$$\sum_i(\theta_i^{(1)})^2=2\sum_i \eta_i^{(2)}+\sum_i\theta_i^{(2)}$$
throughout.
\end{pf}

\subsection{The cohomology of $\ol{M}_{0,0}(\pp^\infty,3)$}

Introduce coordinates $(q_i)$ as in the discussion leading up to
Corollary~\ref{tring}. 

Let us denote by $T_n$ the \'etale $\ol{M}_{0,0}(\pp^n,3)$-scheme
constructed in Section~\ref{SecPar}. For $1\leq n<m$ we consider $\pp^n$ as
a subvariety of $\pp^m$ as in Section~\ref{sminf}. We get an induced
commutative diagram
\begin{equation}\label{tntm}
\vcenter{\xymatrix{
T_m\rto & \ol{M}_{0,0}(\pp^m,3)\\
T_n\uto\rto & \ol{M}_{0,0}(\pp^n,3)\uto}}
\end{equation}
The vector field $V$ on $\ol{M}_{0,0}(\pp^m,3)$ restricts to the vector
field $V$ on $\ol{M}_{0,0}(\pp^n,3)$, so there is no ambiguity if we use
the same later $V$ to denote the vector field on these two stacks.

The diagram (\ref{tntm}) induces the commutative diagram of $\cc$-algebras
\begin{equation}\label{ho}
\vcenter{\xymatrix{
\hh^0(T_m,K_V\com)\ar@{->>}[d] &
\hh^0\big(\ol{M}_{0,0}(\pp^m,3),K_V\com\big) \ar@{->}[d]\lto \\
\hh^0(T_n,K_V\com) &
\hh^0\big(\ol{M}_{0,0}(\pp^n,3),K_V\com\big)\lto }}
\end{equation}
In terms of the explicit description given in Corollary~\ref{tring}, the
vertical arrow on the left hand side of (\ref{ho}) is given by the
canonical map
\begin{multline*}
\cc[b,q_1,q_2,q_3,\ell_1,\ell_2,\ell_3]/
(q_1\ell_1,q_2\ell_2,q_3\ell_3,R_{m+1}) \\\longrightarrow
\cc[b,q_1,q_2,q_3,\ell_1,\ell_2,\ell_3]/
(q_1\ell_1,q_2\ell_2,q_3\ell_3,R_{n+1})
\end{multline*}
and is therefore visibly surjective and moreover an isomorphism on the
part of degree less than $n$. The vertical arrow on the right hand side of
(\ref{ho}) is an isomorphism in degrees less than $n$ by
the comparison with $H_{DR}\big(\ol{M}_{0,0}(\pp^m,3)\big)\to
H_{DR}\big(\ol{M}_{0,0}(\pp^n,3)\big)$, which enjoys these properties by
the results of Section~\ref{sminf}.

Just as in Section~\ref{sminf}, we may therefore define limit algebras as
follows:
$$\hh^0(T_\infty,K_V\com)=\bigoplus_p \projectlim_n\hh^0(T_n,K_V\com)_p\,,$$
and
$$\hh^0\big(\ol{M}_{0,0}(\pp^\infty,3),K_V\com\big)=\bigoplus_p
\projectlim_n\hh^0\big(\ol{M}_{0,0}(\pp^n,3),K_V\com\big)_p\,,$$ 
where the subscript $p$ denotes the homogeneous component of degree $p$. As
above, the limits stabilize as soon as $n>p$.  For every $n$ we have a
commutative diagram of $\cc$-algebras
\begin{equation}\label{hoho}
\vcenter{\xymatrix{
\hh^0(T_\infty,K_V\com)\ar@{->>}[d] &
\hh^0\big(\ol{M}_{0,0}(\pp^\infty,3),K_V\com\big) \ar@{->}[d]\lto \\
\hh^0(T_n,K_V\com) &
\hh^0\big(\ol{M}_{0,0}(\pp^n,3),K_V\com\big)\lto }}
\end{equation}

\begin{them}\label{themi}
The canonical morphism
\begin{multline}\label{canmor}
H_{DR}\big(\ol{M}_{0,0}(\pp^\infty,3)\big)=
\hh^0\big(\ol{M}_{0,0}(\pp^\infty,3),K_V\com\big)\\
\longrightarrow \hh^0(T_\infty,K_V\com) 
=\cc[b,q_1,q_2,q_3,\ell_1,\ell_2,\ell_3]/(q_1\ell_1,q_2\ell_2,q_3\ell_3)\,,
\end{multline}
is injective. Moreover,
we have
$$H_{DR}\big(\ol{M}_{0,0}(\pp^\infty,3)\big)=
\cc[b,\sigma_1,\rho,\sigma_2,\tau,\sigma_3]/
\big((\tau^2-\rho\sigma_2),\tau\sigma_3,\rho\sigma_3\big)\,,$$
where, via~(\ref{canmor}), we have the following identifications:
\begin{align*}
\sigma_1& =\ell_1+\ell_2+\ell_3\,,\\
\rho &=q_1^2+q_2^2+q_3^2-2(q_1 q_2+q_1 q_3+q_2 q_3)\,,\\
\sigma_2 &=\ell_1^2+\ell_2^2+\ell_3^2-2(\ell_1 \ell_2+\ell_1 \ell_3+\ell_2
\ell_3)\,,\\ 
\tau&=\ell_1(q_3-q_2) + \ell_2(q_1-q_3) +\ell_3(q_2-q_1)\,,\\
\sigma_3 &= \ell_1\ell_3\ell_3\,.
\end{align*}
Note that $\rho=\sum\eta_i^{(2)}+\sum\theta_i^{(2)}$ and
$\tau=\sum\ell_i\theta_i^{(1)}$ are the functions that appeared in
Corollary~\ref{second}.
\end{them}
\begin{pf}
Let $A$ be the graded $\cc$-subalgebra of
$\hh^0\big(\ol{M}_{0,0}(\pp^\infty,3),K_V\com\big)$ 
generated by the Chern classes $c_1(E_1)$, $c_1(E_2)$, $c_2(E_1)$,
$c_2(E_2)$, $c_2(E_3)$ and $c_3(E_1)$. Let $B$ the the image of $A$ in
$\hh^0(T_\infty,K_V\com)$. Then using Corollaries~\ref{first}
and~\ref{second} and a direct calculation of $c_3(E_1)$ one shows that $B$
is generated by 
$$b,\sigma_1,\rho,\sigma_2,\tau,\sigma_3\in
\cc[b,q_1,q_2,q_3,\ell_1,\ell_2,\ell_3]/(q_1\ell_1,q_2\ell_2,q_3\ell_3)\,.$$
Now in \cite{getzler} it is proved that the limit
$$\lim_{n\to\infty}\dim H_{DR}^p\big(\ol{M}_{0,0}(\pp^n,3)\big)$$
exists, and is equal to the coefficient of $t^p$ in~(\ref{ezra}).
So from the purely algebraic Lemma~\ref{alglem}, we conclude that
\begin{align*}
\dim H^p_{DR}\big(\ol{M}_{0,0}(\pp^\infty,3)\big) & \geq \dim A_p\\
&\geq \dim B_p\\
&=\lim_{n\to\infty} \dim H_{DR}^p\big(\ol{M}_{0,0}(\pp^n,3)\big)\,.
\end{align*}
Hence we must have equality throughout, which proves the theorem.
\end{pf}

\begin{lem}\label{alglem}
The homomorphism of graded $\cc$-algebras
\begin{multline*}
\cc[b,\sigma_1,\rho,\sigma_2,\tau,\sigma_3]/
\big((\tau^2-\rho\sigma_2),\tau\sigma_3,\rho\sigma_3\big)\\
\longrightarrow
\cc[b,q_1,q_2,q_3,\ell_1,\ell_2,\ell_3]/(q_1\ell_1,q_2\ell_2,q_3\ell_3)
\end{multline*}
defined by the formulas of Theorem~\ref{themi} is injective. Moreover, the
Hilbert series of $\cc[b,\sigma_1,\rho,\sigma_2,\tau,\sigma_3]/
\big((\tau^2-\rho\sigma_2),\tau\sigma_3,\rho\sigma_3\big)$ is given by
\begin{equation}\label{ezra}
\frac{1+t+2t^2+2t^3+2t^4}{(1-t)(1-t^2)^2(1-t^3)}\,.
\end{equation}
\end{lem}
\begin{pf}
Even a computer can prove this lemma, using {\em Macaulay 2}, for
example. See~\cite{M2}.
\end{pf}

\begin{rmk}
Using the results of Section~\ref{chern} it is possible to express the
characteristic classes of the bundles $E_m$ in terms of the generators given
in Theorem~\ref{themi}.  For the first two Chern classes we get
$$c_1(E_m)=
m (3m+1)\,b+\frac{1}{2}m (m+1 )\,\sigma_1\,,$$
\begin{multline*}
c_2(E_m)=\\
\frac{3}{2}{m}^{3} (3m+1 )\,{b}^{2}
+\frac{3}{2}{m}^{3} (m+1)\,b\sigma_{{1}}
+\frac{1}{24}m (m+1 ) (3{m}^{2}+m-1 )\,\sigma_1^{2}\\
-\frac{1}{2}m ( 3m+1 ) (3m+2 )\,\rho 
- \frac{1}{24}m (m+1 ) (2m+1)\,\sigma_2
-\frac{1}{6}m (m+1 ) (7m+2 )\,\tau\,.
\end{multline*}
\end{rmk}

\begin{cor}
For every finite value of $n$ there exists a natural morphism of
$\cc$-algebras
\begin{equation}\label{morex}
\cc[b,\sigma_1,\rho,\sigma_2,\tau,\sigma_3]/
\big((\tau^2-\rho\sigma_2),\tau\sigma_3,\rho\sigma_3\big)
\longrightarrow 
H_{DR}\big(\ol{M}_{0,0}(\pp^n,3)\big)\,,
\end{equation}
which is an isomorphism in degrees less than $n$. 
\end{cor}

\begin{cor}\label{generation}
For every $n$, the cohomology ring $H_{DR}\big(\ol{M}_{0,0}(\pp^n,3)\big)$ 
is generated by the 
characteristic classes of the bundles $E_m$, $m\geq1$ (In fact the first
three $E_m$ and the first three Chern classes suffice).
\end{cor}
\begin{pf}
Here we use Remark~\ref{remarksurjective}, from which we know that
(\ref{morex}) is an algebra epimorphism.
\end{pf}

\subsection{A conjectural presentation for finite $n$}

Let $A$ be the $5\times5$-matrix
$$A=
\left (\begin {array}{ccccc} b+\sigma_{{1}}&0&9\,\rho+\frac{1}{4}(
\sigma_{{2}}-{\sigma_{{1}}}^{2})+3\,\tau&0&\sigma_{{3}}
\\\noalign{\medskip}0&b+\frac{1}{2}\,\sigma_{{1}}&4\,\rho+\frac{1}{2}\,\tau&
\rho&-\frac{1}{2}\,\sigma_{{1}}\rho\\\noalign{\medskip}1&0&b&0&0
\\\noalign{\medskip}0&1&0&b&0\\\noalign{\medskip}0&0&1&1&b
\end {array}\right )$$
and $G_1$ the column vector
$$G_1=
\begin{pmatrix}
b\,\sigma_{{1}}+18\,\rho+\frac{1}{2}(\sigma_{{2}}+\sigma_1^2)+6\,\tau\\
9\,\rho+\tau\\ 2\,b+\sigma_{{1}}\\ b\\ 3
\end{pmatrix}\,.$$

\begin{conjecture}\label{maybe}
We conjecture that
$$H_{DR}\big(\ol{M}_{0,0}(\pp^n,3)\big)=
\cc[b,\sigma_1,\rho,\sigma_2,\tau,\sigma_3]/
\big((\tau^2-\rho\sigma_2),\tau\sigma_3,\rho\sigma_3,G_{n+1}\big)\,,$$
where $G_{n+1}=A^n G_1$.  Note that the degrees of the moving
relations given by $G_{n+1}$ are $n$, $n+1$, $n+1$, $n+2$ and $n+2$.
\end{conjecture}

\begin{lem}
There exists a morphism of $\cc$-algebras
\begin{multline}\label{contilde}
\cc[b,\sigma_1,\rho,\sigma_2,\tau,\sigma_3]/
\big((\tau^2-\rho\sigma_2),\tau\sigma_3,\rho\sigma_3,G_{n+1}\big)\\
\longrightarrow
\cc[b,q_1,q_2,q_3,\ell_1,\ell_2,\ell_3]/
(q_1\ell_1,q_2\ell_2,q_3\ell_3,R_{n+1})
\end{multline}
defined by the formulas of Theorem~\ref{themi}.
\end{lem}
\begin{pf}
Consider the $5\times4$-matrix $H$ given by
$$\left(\begin{array}{cccc}
\sigma_1 &
\rho+\ell_1^2-q_1\lambda_1+5\ell_1\theta_1 &
\rho+\ell_2^2-q_2\lambda_2+5\ell_1\theta_2 &
\rho+\ell_3^2-q_3\lambda_3+5\ell_1\theta_3 \\
0 &
3\rho + \ell_1\theta_1 &
3\rho + \ell_2\theta_2 &
3\rho + \ell_3\theta_1 \\
2 &
\ell_1 + 2\theta_1 &
\ell_2 + 2\theta_2 &
\ell_3 + 2\theta_3 \\
1 & \theta_1 & \theta_2 & \theta_3 \\
0 & 1 & 1 & 1
\end{array}\right)$$
where $\lambda_1=\ell_3-\ell_2$, $\lambda_2=\ell_1-\ell_3$,
$\lambda_3=\ell_2-\ell_1$ and $\theta_1=q_3-q_2$, $\theta_2=q_1-q_3$,
$\theta_3=q_2-q_1$.
This matrix satisfies 
\begin{align*}
AH \equiv HE &\bmod (q_1\ell_1,q_2\ell_2,q_3\ell_3) \\
HR_1  \equiv G_1&\bmod (q_1\ell_1,q_2\ell_2,q_3\ell_3)\,.
\end{align*}
By induction, this implies 
$$G_{n+1}\equiv HR_{n+1}\bmod (q_1\ell_1,q_2\ell_2,q_3\ell_3)\,,$$
which implies the lemma.
\end{pf}

We have the following evidence for our conjecture:

\begin{prop}\label{evidence}
Conjecture~\ref{maybe} holds for $n\leq 5$.
\end{prop}
\begin{pf}
The following two (purely algebraic) conjectures imply
Conjecture~\ref{maybe}. 

Conjecture 1. The Hilbert series of the graded ring
$$C_n=\cc[b,\sigma_1,\rho,\sigma_2,\tau,\sigma_3]/
\big((\tau^2-\rho\sigma_2),\tau\sigma_3,\rho\sigma_3,G_{n+1}\big)$$
is given by
\begin{align}\label{ezra2}
\left.\frac{(1-t^n)(1-t^{n+1})}{(1-t)(1-t^2)^2(1-t^3)}
\right[&t^{2n+3}(2+2t+2t^2+t^3+t^4)\\
-&t^{n+1}(1+3t+4t^2+4t^3+3t^4+t^5)\nonumber\\
+&(1+t+2t^2+2t^3+2t^4)\left]\phantom{\frac{1}{1}}\right..\nonumber
\end{align}

Conjecture 2. Let $B_n$ be the image of the ring morhpism
(\ref{contilde}). Then for all $p\leq n+2$ the dimension of the graded
piece of degree $p$ of $B_n$ is given by the  coefficient of $t^p$ in
Formula~(\ref{ezra2}). 

Both conjectures can be verified algorithmically for any given value
of $n$.  The authors did this using {\em Macaulay 2} (see \cite{M2})
for $n\leq5$.  In fact, the calculations indicate that $p\leq n+2$ is
not the best possible bound.

It is proved in \cite{getzler}, that (\ref{ezra}) is the Poincar\'e
polynomial of $\ol{M}_{0,0}(\pp^n,3)$. 

Thus, using Conjecture~2, the morphism of
$\cc$-algebras
$$H_{DR}\big(\ol{M}_{0,0}(\pp^n,3)\big)\longrightarrow B_n$$ 
is an isomorphism in degrees $p\leq n+2$. Consider the epimorphism
(\ref{contilde}) $C_n\to B_n$. By Conjecture~1, it is also an
isomorphism in degrees $p\leq n+2$.
$$\xymatrix{
& C_n \ar@{..>}[dl] \ar@{->>}[d] \\
H_{DR}\big(\ol{M}_{0,0}(\pp^n,3)\big) \rto & B_n}
$$
Since all relations in $C_n$ are of degrees $p\leq n+2$, we see that
$C_n\to B_n$ lifts as indicated in the diagram. Now that we have an
algebra morphism $C_n\to H_{DR}\big(\ol{M}_{0,0}(\pp^n,3)\big)$, we
note that it is surjective by Remark~\ref{remarksurjective} and hence
an isomorphism, because both rings have the same Hilbert series.
\end{pf}

\subsection{The degree 2 case}

Recall that the coordinates on
$T=\aaa^2\times(\aaa^1\times\aaa^2)^{n-1}$ are called
$(b,q,\ldots,b_\nu,r_{\nu1},r_{\nu2},\ldots)$.  Accordingly, the
vector field $V$ has components 
$$V=V_b\frac{\del}{\del b}
+V_q\frac{\del}{\del q}
+\sum_{\nu=2}^n\begin{pmatrix} \frac{\del}{\del b_\nu}
\frac{\del}{\del r_{\nu1}}\frac{\del}{\del r_{\nu2}}\end{pmatrix}
\cdot V_\nu\,,$$
where $V_b$ and $V_q$ are regular functions on $T$, and $V_\nu$, for
$\nu=2,\ldots,n$ is a column vector thereof.

\begin{prop}
With this notation, the vector field $V$ on $T$ is given by 
\begin{align*}
V_b &= b_2-b^2-4q\,,\\
V_q &= q(r_{21}+r_{22}-4b)\,,\\
V_\nu &= r_{\nu+1}-E r_\nu\,,\quad\text{for $\nu=2,\ldots,n-1$}\\
V_n &=-Er_n\,,
\end{align*}
where for $\nu=2,\ldots,n$ 
$$r_\nu=\begin{pmatrix}b_{\nu}\\r_{\nu1}\\r_{\nu2}\end{pmatrix}$$
and
$$E=\begin{pmatrix}b &2q &2q \\ 1&
r_{21}-b&0\\1&0&r_{22}-b\end{pmatrix}\,.$$
\end{prop}
\begin{pf}
Let $T'\subset T$, be the locus over which we have not blown up, so
that the universal curve over $T'$ is equal to $T'\times\pp^1$.
It suffices to exhibit a vertical vector field
$\ol{U}\in\Gamma(T'\times\pp^1,\tT_{\pp^1})$, such that, if we define
$V$ by the formulas of the proposition, then $Df(V+\ol{U})=W$, where
$W$ is the vector field on $\pp^n$ described in Example~\ref{stupex}.
This means that 
$$V(\varphi_\nu)+\ol{U}(\varphi_\nu)=\varphi_{\nu+1}
-\varphi_1\varphi_\nu\,,$$  
for $\nu<n$ and 
$$V(\varphi_n)+\ol{U}(\varphi_n)=-\varphi_1\varphi_n\,.$$
Here $\varphi_1=b+\frac{1}{s}+qs$ and
$\varphi_\nu=b_\nu+r_{\nu1}\frac{1}{s} +r_{\nu2}qs$, for $\nu>1$. One
checks that 
$$\ol{U}=\frac{\del}{\del s}+(2b-r_{21})s\frac{\del}{\del s}
-qs^2\frac{\del}{\del s}$$
solves this problem.
\end{pf}

Let us write $R_1=\begin{pmatrix}b\\1\\1\end{pmatrix}$. 

\begin{cor}
With this abbreviation, we have
$$\hh^0(T,K_V\com) =
\cc[b,q,r_1,r_2]/\big(q(r_1+r_2-4b),R_{n+1}\big)\,,$$
where $R_{n+1}=E^n R_1$ and we have dropped the first index of the
$r_{2i}$. 

The degrees of $b$, $r_1$, $r_2$ are 1, the degree of $q$ is 2 and the
degrees of the three components of $R_{n+1}$ are $n+1$, $n$ and
$n$. The symmetric group $S_2$ acts by switching $r_1$ and $r_2$ and
switching the last two components of $R_{n+1}$.
\end{cor}

Now we determine the Chern classes of the vector bundles $E_m=\pi\lst
f\upst\O(m)$. 

\begin{prop}\label{deg2chern}
We have
$$c_1(E_m)=-mb+\frac{1}{2}m(m+1)(r_1+r_2)\,,$$
and
$$c_2(E_m)=\frac{1}{2}\big(c_1(E_m)^2-\tr M_m^2\big)\,,$$
where
\begin{multline*}
\tr M_m^2=\frac{2}{3}m(m^2+6m+2) b^2 -\frac{1}{3}m(m+1)(m+2)
b(r_1+r_2)\\
+\frac{1}{6}m(m+1)(2m+1)(8q+r_1^2+r_2^2)\,.
\end{multline*}
\end{prop}
\begin{pf}
We can trivialize $E_m$ over $T$ with the basis
$(1,f_1,f_2,\ldots,f_1^m,f_2^m)$, where
$f_1\in\Gamma\big(\tilde{C}_T,\O(\tilde{x}_1)\big)$ and 
$f_2\in\Gamma\big(\tilde{C}_T,\O(\tilde{x}_2)\big)$ are the
meromorphic functions $f_1=\frac{1}{s}$, and $f_2=qs$. Then we have
for the action of $V$ on this basis:
\begin{align*}
\tilde{V}(1)&=m(b+f_1+f_2)\,,\\
\tilde{V}(f_1^\mu)&=(\mu+m)qf_1^{\mu-1}+((m-2\mu)b+\mu
r_1)f_1^\mu+(m-\mu) f_1^{\mu+1}\,,\\
\tilde{V}(f_2^\mu)&=(\mu+m)qf_2^{\mu-1}+((m-2\mu)b+\mu r_2) f_2^\mu +
(m-\mu)f_2^{\mu+1}\,.
\end{align*}
This leads directly to the formulas of the proposition.
\end{pf}

We are now ready to compute cohomology rings. Let us begin with the
case $n=\infty$. Denote by $\hh^0(T_\infty,K_V\com)$ the limit of
$\hh^0(T,K_V\com)$ as $n$ goes to $\infty$. Thus
$$\hh^0(T_\infty,K_V\com)=\cc[b,q,r_1,r_2]/\big(q(r_1+r_2-4b)\big)\,.$$

\begin{prop}\label{deg2infty}
The canonical morphism 
\begin{equation}\label{deg2}
\hh^0\big(\ol{M}_{0,0}(\pp^\infty,2),K_V\com\big)\longrightarrow
\hh^0(T_\infty,K_V\com)
\end{equation}
is injective. We have
$$H_{DR}\big(\ol{M}_{0,0}(\pp^\infty,2)\big)=\cc[b,t,k]\,,$$
where, under (\ref{deg2}), we have
$b\mapsto b$, $t\mapsto r_1+r_2-2b$ and $k\mapsto 4q-(b-r_1)(b-r_2)$.
\end{prop}
\begin{pf}
This is proved the same way as Theorem~\ref{themi}. One uses that the
Betti numbers of $\ol{M}_{0,0}(\pp^n,2)$ stabilize to the coefficients
of 
$$\frac{1}{(1-t)^2(1-t^2)}$$
(see \cite{getzler}).
\end{pf}

\begin{cor}\label{td1}
For every finite $n$, the cohomology ring of $\ol{M}_{0,0}(\pp^n,2)$
is generated by the Chern classes of the bundles $E_m$. In fact,
$c_1(E_1)$, $c_1(E_2)$ and $c_2(E_1)$ will suffice.
\end{cor}
\begin{pf}
It follows from Remark~\ref{remarksurjective} that the cohomology ring
of $\ol{M}_{0,0}(\pp^n,2)$ is generated by $b$, $t$ and $k$, which
translates easily into the claim about Chern classes, by virtue of the
formulas of Proposition~\ref{deg2chern}. 
\end{pf}

Now we are ready to determine the cohomology ring of
$\ol{M}_{0,0}(\pp^n,2)$ for finite $n$.

\begin{prop}\label{td2}
We have
$$H_{DR}\big(\ol{M}_{0,0}(\pp^n,2)\big)=\cc[b,t,k]/(G_{n+1})\,.$$
Here $G_{n+1}$ stands for three relations in degrees $n$, $n+1$ and
$n+2$, defined recursively by the matrix equation $G_{n+1}=A^n G_1$,
where
$$A=\begin{pmatrix} b & 0 & 0\\1&0&k\\0&1&t\end{pmatrix}$$
and 
$$G_1=\begin{pmatrix} b(2b-t)\\2b-t\\2\end{pmatrix}\,.$$
\end{prop}
\begin{pf}
The Poincar\'e polynomial of $\ol{M}_{0,0}(\pp^n,2)$, which we glean
from \cite{getzler} is
\begin{equation}\label{deg2getzler}
\frac{(1-t^n)(1-t^{n+1})(1-t^{n+2})}{(t-1)^2(1-t^2)}\,.
\end{equation}
Now the proof proceeds in the same way as the proof of
Proposition~\ref{evidence}, except for that we can prove the two
requisite facts. They follows as lemmas.
\end{pf}

In the following, we will abbreviate the three components of $G_{n+1}$
as $u_{n+1}$, $v_{n+1}$ and $w_{n+1}$. In particular,
$u_{n+1}=b^{n+1}(2b-t)$. Similarly, the three components of $R_{n+1}$
are $x_{n+1}$, $y_{n+1}$ and $z_{n+1}$.

\begin{lem}
The Hilbert series of the graded ring
$\cc[b,t,k]/(G_{n+1})$
is given by~(\ref{deg2getzler}).
\end{lem}
\begin{pf}
First note that the only common zero of the three components of $G_{n+1}$
is $0$. One way to prove this follows. Assume that $u_{n+1}=0$,
$v_{n+1}=0$ and 
$w_{n+1}=0$. By the explicit description $b^{\mu}(2b-t)$ of
$u_{\mu}$ we see that we get $u_\mu=0$, for all $\mu\geq1$. Then we
have $v_{\mu+1}=kw_{\mu}$, for all $\mu$, and the assumption that
$k\not=0$, forces $w_n=0$. The conjunction of $w_{n+1}=0$ and $w_n=0$
forces $v_n=0$. Proceeding inductively, we see that $w_1=0$, which
contradicts the fact that $w_1=2$. Thus we can conclude that
$k=0$. Then all $v_\nu=0$ and we have $w_{\mu+1}=tw_\mu$, and so
$w_\mu=2t^{\mu-1}$. This implies that $t=0$. Finally, $u_{n+1}=0$
then implies that $b=0$.

We conclude that  $\dim\cc[b,t,k]/(G_{n+1})=0$, which implies that
$G_{n+1}$ is a regular sequence.  Thus we can use the Koszul complex
of $G_{n+1}$ to compute the Hilbert series of this quotient ring. We
get~(\ref{deg2getzler}).
\end{pf}

\begin{lem}
The formulas of Proposition~\ref{deg2infty} define a $\cc$-algebra
morphism
$$\cc[b,t,k]/(G_{n+1})\longrightarrow
\cc[b,q,r_1,r_2]/\big(q(r_1+r_2-4b),R_{n+1}\big)\,.$$
This morphism is injective in degrees $\leq n+2$.
\end{lem}
\begin{pf}
The matrix
$$H=\begin{pmatrix}4b-r_1-r_2 &0&0\\2& b-r_2 &
b-r_1\\0&1&1\end{pmatrix}$$
satisfies $HE\equiv AH\bmod q(r_1+r_2-4b)$ and $H
R_1=G_1$.
By induction, we conclude that $HR_{n+1}\equiv
G_{n+1}\bmod q(r_1+r_2-4b)$.  This proves the existence of the
morphism.
It is injective in degrees less than $n$, because
$\cc[b,t,k]\to\cc[b,q,r_1,r_2]/q(r_1+r_2-4b)$ is injective. It remains
to prove injectivity in degrees $n$, $n+1$ and $n+2$.

The invariant subring will be helpful. It is given by
\begin{multline*}
\left(\cc[b,q,r_1,r_2]/\big(q(r_1+r_2-4b),R_{n+1}\big)\right)^{S_2}\\
=\cc[b,t,k,q]/\big(q(2b-t),x_{n+1},v_{n+1},w_{n+1}\,.
\end{multline*}
The proof that this is indeed the ring of invariants proceeds as
follows. Fist note that the invariant subring of $\cc[b,q,r_1,r_2]$ is
the polynomial ring $\cc[b,t,k,q]$. Then, using the abbreviation
$B=\cc[b,q,r_1,r_2]$, we consider the exact sequence
$$B^4\longrightarrow B\longrightarrow
B/q(r_1+r_2-4b),x_{n+1},y_{n+1},z_{n+1}\longrightarrow 0$$
of $B$-modules. Taking invariants of this sequence, gives us an exact
sequence of $B^{S_2}$-modules. Since the invariants of $B^4$ are
generated by $(1,0,0,0)$, $(0,1,0,0)$, $(0,0,1,1)$ and
$(r_1-r_2)(0,0,1,-1)$, we conclude that the invariant subring of
$B/q(r_1+r_2-4b),x_{n+1},y_{n+1},z_{n+1}$ is equal to
$\cc[b,t,k,q]/q(2b-t),x_{n+1},(r_1-r_2)(y_{n+1}-z_{n+1}),(y_{n+1}+z_{n+1})$.
That this is equivalent to the above presentation, follows from
$HR_{n+1}\equiv G_{n+1}$.

So we are now reduced to proving that 
\begin{multline*}
\cc[b,t,k]/\big(b^{n+2}(2b-t),v_{n+1},w_{n+1}\big)\\
\longrightarrow\cc[b,t,k,q]/\big(q(2b-t),x_{n+1},v_{n+1},w_{n+1}\big)
\end{multline*}
is injective in degrees $n$, $n+1$ and $n+2$. 
In the presentation
$\cc[b,t,k,q]/\big(q(2b-t),x_{n+1},v_{n+1},w_{n+1}\big)$, we may
replace $x_{n+1}$ by $b^{n+1}+q\tilde{x}_{n+1}$, where
$\tilde{x}_{n+1}\in\cc[b,t,k]$ and is recursively defined by
$\tilde{x}_1=0$, $\tilde{x}_{n+1}=b\tilde{x}_n+2w_n$.
Once we have done this, the proof proceeds as follows. 

Let $a\in\cc[b,t,k]$ be an element mapping to zero in
$\cc[b,t,k,q]/\big(q(2b-t),b^{n+1}+q\tilde{x}_{n+1},v_{n+1},w_{n+1}$.
Thus
\begin{equation}\label{expandq}
a=e(q)q(2b-t)+f(q)(b^{n+1} +q\tilde{x}_{n+1}) + g(q)v_{n+1} +
h(q)w_{n+1}\,.
\end{equation}
We may assume that every term in this equation is homogeneous. 
Let us expand in powers of $q$.  The constant term gives us
$$a=f_0b^{n+1} +g_0 v_{n+1} + h_0 w_{n+1}\,.$$
If $\deg a=n$, then for degree reasons it follows that $f_0=0$,
$g_0=0$ and $h_0\in\cc$. Thus $a\in(w_{n+1})$, and we are done.

If $\deg a =n+1$, then for degree reasons $f_0\in\cc$. Let us
consider the linear term of~(\ref{expandq}). 
$$0=e_0(2b-t) + f_0\tilde{x}_{n+1}\,.$$
Since $(2b-t)$ does not divide $\tilde{x}_{n+1}$ (see below), this
implies that $f_0=0$. Thus $a\in(v_{n+1},w_{n+1})$, and we are done.

Finally, let $\deg a=n+2$. Then $\deg f_0=1$. This time the linear
term of~(\ref{expandq}) is
$$0= e_0(2b-t) + f_0\tilde{x}_{n+1} + h_1w_{n+1}\,,$$
where $h_1\in\cc$.
Since $\tilde{x}_{n+1}$ does not divides $w_{n+1}$ modulo $2b-t$, this
implies that $h_1=0$. This, in turn, implies that $(2b-t)$ divides
$f_0$. Hence $a\in(b^{n+1}(2b-t),v_{n+1},w_{n+1})$, and we are done.

It remains to check our two claims to the effect that if $t=2b$, then
$\tilde{x}_{n+1}\not=0$ and does not divide $w_{n+1}$.  These can
easily be proved by further setting $k=0$, to solve the recursions.
\end{pf}

\begin{numrmk}\label{final}
It might be instructive to consider the case $d=1$ in this context.
Recall
that  $\ol{M}_{0,0}(\pp^n,1)$ is the Grassmannian of lines in
$\pp^n$. We let $T=\aaa^{2n-2}$, with coordinates
$(b_2,r_2,\ldots,b_n,r_n)$. The universal curve over $T$ is $T\times
\pp^1$.  Let $s$ be an affine coordinate on $\pp^1$, then the
universal map $f:T\times\pp^1\to\pp^n$ is given by 
$$s\longmapsto\langle1,s,b_2+r_2s,\ldots,b_n+r_ns\rangle\,.$$
The auxiliary vertical vector field $\ol{U}$ is
$(b_2+r_2s-s^2)\frac{\del}{\del s}$. This leads to the vector field
$V$ on $T$ in the following form:
\begin{align*}
V_\nu&=\begin{pmatrix} b_{\nu+1}\\r_{\nu+1}\end{pmatrix}-
\begin{pmatrix} 0 & b_2\\1& r_2\end{pmatrix} \begin{pmatrix} b_\nu\\
r_\nu\end{pmatrix}\,,\quad\text{ for $\nu=2,\ldots,n-1$}\,\\
V_n&=-
\begin{pmatrix} 0 & b_2\\1& r_2\end{pmatrix} \begin{pmatrix} b_n\\
r_n\end{pmatrix}\,.
\end{align*}
Here $V_\nu$ is the column vector of functions on $T$ defined by 
$$V=\sum_{\nu=2}^{n} \begin{pmatrix}\frac{\del}{\del b_\nu} &
\frac{\del}{\del r_\nu} \end{pmatrix} V_\nu\,.$$
The only zero of the vector field $V$ is the origin of $T$. Thus we
have
\begin{multline*}
H^\ast(\ol{M}_{0,0}(\pp^n,1),\cc)
=\hh^0\big(\ol{M}_{0,0}(\pp^n,1),K_V\com\big)\\
= \Gamma (T,\O)/V\Gamma(T,\Omega) = \cc[b_2,r_2]/R_{n+1}\,,
\end{multline*}
where 
$$R_{n+1}=\begin{pmatrix} 0 & b_2\\1& r_2\end{pmatrix}^n
\begin{pmatrix} 0\\1\end{pmatrix}\,.$$
Note that the degree of $b_2$ is
2 and the degree of $r_2$ is 1.  To find the geometric significance of
$r_2$ and $b_2$, we calculate the action of $V$ on the vector bundles
$E_m={p_T}\lst f\upst\O(m)$ over $T$.  A basis is provided by
$(1,s,\ldots,s^m)$ and $V$ acts by $\tilde{V}(s^\mu)=\mu b_2
s^{\mu-1}+\mu r_2 s^\mu + (m-\mu) s^{\mu+1}$.  This leads to
$$c_1(E_m)=\frac{1}{2}m(m+1) r_2\,.$$
Moreover, $c_1(E_1)=r_2$ and $c_2(E_1)=-b_2$. 

In this case, there is a nice interpretation of the relations
$R_{n+1}$ in the cohomology ring. On $\ol{M}_{0,0}(\pp^n,1)$ there is
an exact sequence of vector bundles
$$\ses{K}{}{\O^{n+1}}{}{E_1}\,,$$
simply obtained by pushing forward and pulling back the tautological
sequence on $\pp^n$.  The relations $R_{n+1}$ are equivalent to the
equation
$$c_t(K)c_t(E_1)=1$$
among total Chern classes, which follows from the above short exact
sequence.  This short exact sequence does not generalize to the cases
of higher degree, and we know no such simple motivation of our
relations.
\end{numrmk}


\end{document}